\pdfoutput=1
\documentclass[%
 aip, amsmath,amssymb, reprint,twocolumn]{revtex4-1}
\usepackage{graphicx}
\usepackage{dcolumn}
\usepackage{bm}
\usepackage{amsmath}
\usepackage{amsfonts}
\usepackage{amssymb}
\usepackage{cancel}
\usepackage{graphicx}
\usepackage{float}
\usepackage{wrapfig}
\usepackage{caption}
\usepackage{subcaption}
\usepackage{mathtools}

\usepackage{color}
\usepackage{color}
\newcommand{\chg}[1]{#1}
\newcommand{\R}{\mathbb{R}}

\newcommand{\C}{\mathbb{C}}
\renewcommand{\d}{\mathrm{d}}
\newcommand{\e}{\mathrm{e}}
\renewcommand{\i}{\mathrm{i}}
\newcommand{\tmax}{\tau_{\max}}

\newcommand{\re}{\operatorname{Re}}
\newcommand{\im}{\operatorname{Im}}
\newcommand{\id}{I}

\newcommand{\ind}{\operatorname{ind}}
\newcommand{\lspan}{\operatorname{span}}

\newtheorem{theorem}{Theorem}[section]
\newtheorem{lemma}[theorem]{Lemma}

\newtheorem{proposition}[theorem]{Proposition}

\begin{document}


\title{Local bifurcations in differential equations with
  state-dependent delay} 



\author{Jan Sieber}
\email[]{J.Sieber@exeter.ac.uk}
\affiliation{University
  of Exeter, EPSRC Centre for Predictive Modelling in Healthcare,
  University of Exeter, Exeter, EX4 4QJ, UK}


\date{\today}

\begin{abstract}
  A common task when analysing dynamical systems is the determination
  of normal forms near local bifurcations of equilibria. As most of
  these normal forms have been classified and analysed, finding which
  particular class of normal form one encounters in a numerical
  bifurcation study guides follow-up computations.

  This paper builds on normal form algorithms for equilibria of delay
  differential equations with constant delay that were developed and
  implemented in DDE-Biftool recently. We show how one can extend
  these methods to delay-differential equations with state-dependent
  delay (sd-DDEs). Since higher degrees of regularity of local center
  manifolds are still open for sd-DDEs, we give an independent (still
  only partial) argument which phenomena from the truncated normal
  must persist in the full sd-DDE. In particular, we show that all
  invariant manifolds with a sufficient degree of normal hyperbolicity
  predicted by the normal form exist also in the full sd-DDE.
\end{abstract}


\maketitle 

\begin{quotation}
  Delay-differential equations (DDEs) arise frequently in models where
  the evolution of the system depends also on its values in the
  past. Typical examples arise in control (delays in feedback loops),
  optics (delayed feedback effects from external light reflections),
  mechanical engineering (effects from previous rotations in turning
  processes), or Earth sciences (El Ni{\~n}o caused by delayed
  feedback from waves across oceans).

  The typical approach to studying DDEs is to consider them as a
  dynamical systems for which the state is a history segment
  (in our case on a bounded history interval). Several mathematical
  problems occur when the length of the delay depends on the state of
  the system, called sd-DDEs. In this case the state of the dynamical
  system at time $t$ does not depend smoothly on its initial
  condition. This makes many of the standard tools of dynamical
  systems theory inapplicable at first sight. In particular normal
  form theory requires expansion of the right-hand side to higher
  orders.

  This paper demonstrates that normal forms can still be computed for
  a general class of sd-DDEs with discrete delays. We show that the
  computational procedure developed by Janssens, Wage, Bosschaert and
  Kuznetsov\cite{J10,W14,B16,BJK1x} for DDEs with constant delays can
  be generalized to sd-DDEs. We also give a justification for the
  computed normal forms, explaining why all normally hyperbolic
  manifolds present in the normal form also appear in the full
  sd-DDE. The justification is based on an approach recently taken by
  Humphries \emph{et al}\cite{HCK16} in a numerical bifurcation study of a
  prototypical sd-DDE.
\end{quotation}

\section{Introduction}
Delay-differential equations (DDEs) are a class of differential
equations where the derivative at the current time $t$ may depend on
any value of the state in the past. \chg{This paper focusses on those
  case where the dependence is on states from a limited time interval
  $[t-\tmax,t]$ in the past.} They are a particularly common and
well-studied subclass of so-called functional-differential
equations\cite{HL93, DGLW95}. Mathematically, DDEs are dynamical
systems with an infinite-dimensional phase space, since the
appropriate initial value is a prescribed piece of history of the
physical variable on an interval $[-\tmax,0]$. A typical choice of
phase space is the space of $n$-dimensional continuous functions on
$[-\tmax,0]$, written as $C^0([-\tmax,0];\R^n)$ with the maximum norm
(short $C^0$). The right-hand side is given by a functional
$F:C^0\to\R^n$. An example is $F(u)=-u(-\tau)$ for a fixed $\tau>0$
and functions $u$ close to $0$ in $C^0$. Then one will write the
differential equation $\dot u(t)=-u(t-\tau)$ as
\begin{displaymath}
  \dot u(t)=F(u_t)\mbox{,}
\end{displaymath}
where the subscript $t$ indicates a time-shifted history interval. So,
for a function $u:[-\tmax,T]\to\R^n$ and $t\in[0,T]$, $u_t$ is a
function on $[-\tmax,0]$ defined by $u_t(\theta)=u(t+\theta)$.

There is mathematically a large difference between DDEs with constant
delays and DDEs with state-dependent delays. For constant delays, a
framework that poses DDEs as abstract ODE has been developed by Hale
\& Verduyn-Lunel\cite{HL93} and Diekmann \emph{et al}\cite{DGLW95}. In
this framework DDEs of the type $\dot u(t)=F(u_t)$ are smooth
dynamical systems on the phase space $C^0$. That is, the time-$t$ map
$u_0\mapsto u_t$ for fixed $t$, mapping the initial condition $u_0\in
C^0$ to the solution $u_t\in C^0$ at time $t$, is smooth. The
smoothness of the time-$t$ map follows from the smoothness of the
functional $F:C^0\to\R^n$.

This is in contrast to the case when the functional $F$ involves
state-dependent delays. We refer to this type of DDEs as DDEs with
state-dependent delays (short sd-DDEs). \chg{An example is the differential
equation $\dot u(t)=p-u(t+u(t))$ for fixed parameter $p$, for which
the functional $F$ has the form $F:u\mapsto p-u(u(0))$ (for $u$ close
to $p$ and $p<0$). The derivative of the right-hand side $F$ with
respect to its argument $u$ is $\partial F(u)v=-v(u(0))-u'(u(0))v(0)$
if it exists.} Thus, it is undefined for $u\in C^0$ that are not
differentiable. This has the consequence that the standard theory from
textbooks\cite{HL93,DGLW95} for DDEs is not applicable. The currently
most practical statements (for dynamical systems theory) about the
regularity of the time-$t$ map with respect to its initial value are
by Hartung\cite{Hartung11} and Walther\cite{W03}. They are much more
restricted, achieving at best continuous differentiability (once) of
the time-$t$ map. A review by Hartung \emph{et al} from
2006\cite{HKWW06} presents a snapshot of developments regarding
general existence and regularity theory. Section~\ref{sec:review:dde}
summarizes the most relevant results.

\paragraph*{Applications and numerical software}
In parallel to developments in the theory of sd-DDEs, computational
tools have been created to help solving practical problems arising in
engineering and science. The review by Hartung \emph{et
  al}\cite{HKWW06} lists a few classical applications such as control
by echo location\cite{W02}, models for cutting processes with a finite
tool stiffness in directions tangential to the rotating
surface\cite{IBS08,IST07} and the electromagnetic two-body
problem\cite{DGHP10}. Other examples are time-delayed feedback control
where the time-delay is adjusted dynamically\cite{PP11}, and models
for granulopoiesis\cite{CHM16}.

Two common tasks to be performed numerically in applications are
initial-value problem solving (a black-box solver for sd-DDEs including
neutral terms is RADAR5\cite{GH01}) and numerical bifurcation
analysis. Numerical bifurcation analysis tracks branches of equilibria
(constant solutions of $F(u)=0$), periodic orbits (time-periodic
solutions of $\dot u(t)=F(u_t)$) and their bifurcations and linear
stability. Equilibria of sd-DDEs are given by algebraic equations and
periodic boundary-value problems can be reduced to equivalent systems
of smooth algebraic equations\cite{S12}. Thus, numerical computations
of these are feasible in principle and have been implemented in
DDE-Biftool\cite{ELR02,ELS01,ddebiftoolmanual}. Its capabilities for
sd-DDEs with discrete delays (as described in
Section~\ref{sec:discrete}) include:
\begin{itemize}
\item continuation of families of equilibria and computation of their
  stability (present since version 2.0);
\item continuation of codimension-one bifurcations of equilibria (Hopf
  bifurcations and saddle-node bifurcations, present since version
  2.0);
\item continuation of periodic orbits in one parameter and computation
  of their stability (present since version 2.0, completed for the
  class of sd-DDEs with discrete delays described in
  Section~\ref{sec:discrete} in version 3.0);
\item continuation of local codimension-one bifurcations of periodic
  orbits (saddle-node bifurcations, period doubling bifurcations and
  torus bifurcations, present since version 3.0);
\end{itemize}

\paragraph*{Normal forms of local bifurcations}
This paper gives the background on how direct normal form computations
for codimension-one and -two bifurcations of equilibria have been
added for sd-DDEs to the general sd-DDE capabilities. The procedures
are based on the corresponding code and work by Kuznetsov, Janssens,
Wage and Bosschaert\cite{J10,W14,B16,BJK1x} for constant-delay
DDEs. Section~\ref{sec:nf:const} reviews these recent developments for
constant delays. Appendix~\ref{app:locbif} gives more details.

Normal form computations help classify all generic (up to codimension
two) bifurcations into a finite number of well-studied cases. Thus,
they help the systematic numerical exploration in applications. For
example, when a Hopf bifurcation is detected, one may compute the
so-called Lyapunov coefficient which determines to which side the
periodic orbits branch off from the equilibrium (that is, whether the
Hopf bifurcation is \emph{sub-} or \emph{supercritical}, or, using the
terms coined in  engineering, \emph{safe} or
\emph{dangerous}\cite{TS02}). The illustrative example of a linear
position control problem with state-dependent delay in
Section~\ref{sec:pos} shows a typical scenario.

Similarly, when following a Hopf bifurcation in two parameters, one
typically encounters crossings with other Hopf bifurcations (a common
scenario for DDEs). At these so-called Hopf-Hopf interaction points
various branches of secondary bifurcations can be expected depending
on the normal form of the Hopf-Hopf interaction. Humphries \emph{et
  al}\cite{HCK16} studied bifurcations of a scalar sd-DDE in
detail. They encountered several Hopf-Hopf interactions, derived the
normal form on paper, and then followed the predicted secondary
bifurcations, which turned out to exist in the expected directions.

\paragraph*{Justification of normal form expansion in sd-DDEs}
The normal form of most codimension-one and -two bifurcations depends
on expansion terms of order higher than one. Expansion to this degree
is not immediately justifiable for sd-DDEs since the time-$t$ map of
sd-DDEs is only continuously differentiable once. For ordinary
differential equations (ODE), there are precise statements about the
relation between the phase portraits and their bifurcations in
truncated normal forms and the full dynamical system (they depend on
the particular bifurcation)\cite{GH83,K04}. To obtain the same
statements for sd-DDEs one needs that local center manifolds near
equilibria are smooth to the degree required for the expansion terms
in the normal form (for example, to third order for the Hopf
bifurcation). A local center manifold near an equilibrium in a
(sd-)DDE has the form of a graph $h:\R^{n_c}\to C^0([-\tmax,0];\R^n)$.
Here $n_c$ is the number of eigenvalues (counted with multiplicity) of
the linearized DDE on the imaginary axis, and the domain of $h$ is a
coordinate representation of the corresponding eigenspace.  The
smoothness requirement for $h$ refers to two things. First, each
element of the center manifold has to be smooth with respect to its
argument (time), so $h(u_c)\in C^\ell([-\tmax,0];\R^n)$ (the space of
$\ell$ times continuously differentiable functions). Second, the graph
$h$ has to be a smooth map of its argument
$u_c\in\R^{n_c}$. Smoothness of local center manifolds has not been
proven rigorously yet for degrees greater than one. Stumpf
\cite{Stu11} gives a proof of continuous differentiability of
center-unstable manifolds, and shows that it attracts exponentially
all those solutions that stay near the
equilibrium\cite{Stu15a}. However, we prove in
Section~\ref{sec:sd:smooth} that many phenomena predicted by the
normal form must also be present in the sd-DDE. The statement is not
as strong as its classical ODE counterpart such that the availability
of numerical normal form computations provides a motivation to
investigate the smoothness of local center manifolds rigorously.

\section{DDEs with state-dependent delays}
\label{sec:sd-ddes}

\subsection{Discrete state-dependent delays}
\label{sec:discrete}
DDE-Biftool is able to perform bifurcation analysis on a class of
$n$-dimensional systems of delay differential equations with $m-1$
discrete state-dependent delays (sd-DDEs) of the following form:
\begin{align}
  \label{dde:sys_rhs}
  \dot x(t)&=f(x^1,\ldots,x^m,p)\mbox{,\quad where $x^1=x(t)$, and}\\
  \label{dde:sys_tau}
  x^j&=x(t-\tau^j(x^1,\ldots,x^{j-1},p))\mbox{\quad for $j=2,\ldots,m$.}
\end{align}
The integers $n\geq1$ (physical space dimension), $m\geq1$ (number of
delays) and $n_p\geq0$ (number of parameters) are arbitrary. It uses
the convention that $\tau^1=0$ and assumes that the functions
\begin{align}
  \label{dom:sys_rhs}
  f&:\R^{n\times m}\times\R^{n_p}\to\R^n\mbox{,}\\
  \label{dom:tau}
  \tau^j&:\R^{n\times(j-1)}\times\R^{n_p}\to[0,\infty)\mbox{}  
\end{align}
are smooth. The construction \eqref{dde:sys_rhs}--\eqref{dde:sys_tau}
permits arbitrary levels of nesting in the delayed arguments of
$x$. DDE-Biftool does not require an explicit value for the maximal
delay. It computes equilibria and periodic orbits such that the
trajectory $x(t)$ is always compact.

In sections with theoretical considerations we may assume that $n_p=0$
without loss of generality by incorporating the parameters into the
state (appending the equation $\dot p=0$ to \eqref{dde:sys_rhs} and
increasing $n$ to $n+n_p$).

\subsection{General functional differential equations (FDEs) ---
  Review of basic properties}
\label{sec:review:dde}
\paragraph*{Notation and assumptions on the right-hand side}
In the following sections we will use the abbreviation that $C^0$ (or
just $C$) is the space $C([-\tmax,0];\R^n$) of continuous functions on
the interval $[-\tmax,0]$ into $\R^n$ with the norm
\begin{displaymath}
  \|u\|_0=\max\left\{|u(t)|:t\in[-\tmax,0]\right\}\mbox{.}
\end{displaymath}
\chg{Similarly, for any space $D$ of functions on an interval $I\subset\R$
and integer $\ell>0$, we denote the subspace $D^{\ell}$ as the space
of functions which have a $\ell$th derivative in $D$. 
Their
respective norms are
\begin{align*}
  \|u\|_{D^\ell\phantom{,\ell}}&=\max\{\|u\|_D,\|u'\|_D,\ldots,
                                 \|u^{(\ell)}\|_D\}
                      \mbox{.}
\end{align*}}
We also use the phrase, for example, ``$f$ is $C^\ell$\,'' for $f$
being $\ell$ times continuously differentiable in all its arguments.

Basic existence and regularity theory for solutions of sd-DDEs has
been developed for differential equations in the form
\begin{equation}
  \label{eq:genfde}
  \dot u(t)=F(u_t)\mbox{,}
\end{equation}
where $F:C([-\tmax,0];\R^n)\to\R^n$ is a continuous nonlinear
functional \cite{HKWW06}. For a function $u:[-\tmax,T]\to\R^n$ the
notation $u_t$ refers to a time shift of $u$ back to a function on the
interval $[-\tmax,0]$:
\begin{displaymath}
  u_t(\theta)=u(t+\theta)\mbox{\ for $t\in[0,T]$\ 
    and $\theta\in[-\tmax,0]$.}
\end{displaymath}
For the type of equations that
can be treated with DDE-Biftool the functional $F$ (incorporating
parameters into the state variables) has the form
\begin{align}
  \label{fde:F}
  F(u)&=f(u^1,\ldots,u^m)\mbox{,\quad where $u^1=u(0)$, and}\\
  \label{fde:tau}
  u^j&=u(-\tau^j(u^1,\ldots,u^{j-1}))\mbox{\quad for $j=2,\ldots,m$.}
\end{align}
\chg{If the coefficient functions $f$ and $\tau^j$ are $\ell$ times
continuously differentiable, we call such a functional $F$ a
\emph{functional with $C^\ell$ coefficients and $m$ state-dependent
  discrete delays less than $\tmax$}.}

The general conditions on $F$ to ensure existence and regularity of
solutions vary between different papers. A set of conditions that
covers functionals $F$ with discrete state-dependent delays and
$C^\ell$ coefficients and satisfies the assumptions in many
fundamental papers is \emph{mild differentiability}. \chg{Consider a
  continuous functional $F:D\to \R^N$ for some $N\geq1$ and some $D$
  that is a subspace of $C^0(I;\R^N)$ for some interval $I\subset
  \R$. For mild differentiablity of $F$ we require the following two
  conditions.
\begin{enumerate}
\item[(S1)] The functional $F$ is continuously differentiable when
  restricted to the subspace $D^1$. We denote its derivative by
  $\partial F:D^1\to{\cal L}(D^1;\R^N)$.
\item[(S2)] The map
  \begin{displaymath}
    D^1\times D^1\ni(u,v)\mapsto \partial F(u)v\in \R^N
\end{displaymath}
can be extended continuously to the space $D^1\times D$.
\end{enumerate}
We put the argument $v$ of $\partial F$ outside of the bracket to
emphasize that $\partial F$ is linear in $v$. Since $\partial
F:D^1\times D\to\R^N$ is continuous, we can apply the definition
for mild differentiability recursively, treating the pair $(u,v)\in
D^1\times D$ as the single argument of $\partial F$. This leads
naturally to the definition that a functional $F:D\to\R^N$ is
$\ell$ times mildly differentiable if
\begin{enumerate}
\item[(S3)] $\partial F:(D^1\times D)\to\R^N$ is $\ell-1$ times
  mildly differentiable.
\end{enumerate}}

\paragraph*{Scalar illustrative example}
An illustrative example is the sd-DDE \chg{
\begin{align}
  \label{dde:nested}
  \dot x(t)&=p-x(t+x(t))\mbox{,\ that is,\ }\nonumber\\
  \dot u(t)&=F(u_t)\mbox{\ with\ } F(u)=p-u(u(0))\mbox{.}
\end{align}
This corresponds to the choice $f(x,y,p)=p-y$ and $\tau^2(x,p)=-x$} in
\eqref{fde:F}--\eqref{fde:tau} (using letters $x$ and $y$ in the
arguments of $f$ instead of superscripts to avoid confusion with
powers), where we keep \chg{$p=-\pi/2$ fixed for illustration initially.
So, $F$ is a functional with $2$ delays and $C^\infty$
coefficients. The first two derivatives of this functional $F$ are
\begin{align*}
  \partial F(u)v=&-u'(u(0))v(0)-v(u(0))\\
  \partial [\partial F(u,v)](w,z)=&\partial^2F(u)vw+\partial F(u)z\\
  =&-w'(u(0)) v(0)  - v'(u(0)) w(0)\\
  & - u''(u(0)) w(0) v(0)\\
  &- u'(u(0)) z(0) - z(u(0))\mbox{.}
\end{align*}
Note how the second derivative includes differentiation of the first
derivative with respect to $v$ according to our convention such that
it has $4$ arguments (generally, the $\ell$th derivative will have
$2^{\ell+1}$ arguments). We reserve the notation $\partial^jF(u)$ for the
usual $j$-linear form.  The above expressions show that the $\ell$th
derivative of $F$ depends on the lowest $\ell$ derivatives of $u$, on
the lowest $(\ell-1)$ derivatives of the deviation $v$ and $w$, and
only on the values of $z$. So, $\partial^1F$ is continuous in
$C^1\times C^0$ and $\partial[\partial F]$ is continuous in
$(C^2\times C^1\rangle\times(C^1\times C^0)$. Moreover, the map $u
\mapsto \partial F(u,\cdot)$ is continuous as a map, mapping $u\in C^1$ into
the space ${\cal L}(C^1;\R)$ of linear functionals from $C^1$ into
$\R$, but \emph{not} as a map into the space ${\cal L}(C^0;\R)$ of
linear functionals from $C$ into $\R$. The reason for this
discontinuity is the second term $-v(u(0))$: the map
\begin{displaymath}
  [\tmax,0]\ni\theta\mapsto[C^\ell\ni v\mapsto v(\theta)]\in{\cal L}(C^\ell;\R)
\end{displaymath}
is only continuous in $\theta$ if $\ell\geq1$. Mild differentiability
of second order requires that
$(u,v)\mapsto \partial[\partial F(u,v)](\cdot,\cdot)\in{\cal L}(C^2\times
C^1;\R)$
is continuous, which is the case for the right-hand side $F$ in
example \eqref{dde:nested}.}

The example illustrates that the assumptions of mild differentiability
permit dependence of the delays on the state.  \chg{We note that for
  varying $p$, we have to include the equation $\dot p=0$. The
  combined system also satisfies mild differentiability to all
  orders. Equation~\eqref{dde:nested} has an equilibrium at $u=p$,
  which loses its stability in a Hopf bifurcation at $p=-\pi/2$.}  We
will use the above example \eqref{dde:nested} to illustrate various
technical assumptions and difficulties in the following sections.  For
example, the form of the first derivative of $F$ in \eqref{dde:nested}
implies that $F$ is not locally Lipschitz continuous in $C^0$.

\paragraph*{Basic results on solutions of sd-DDEs}
Successive differentiation and application of the chain rule imply
that functionals $F$ with discrete delays and $C^\ell$ coefficients
(in the form of \eqref{fde:F}--\eqref{fde:tau}) satisfy assumptions
(S1--S3) up to the order $\ell$. Thus, all of the following basic
results apply to this class of sd-DDEs with discrete delays.

Walther\chg{\cite{W03,W04}} proved 
that initial value problems (IVPs) have a unique solution $u$ for all
times $t$, or the solution blows up in finite time, if the initial
value $u_0$ lies in the manifold
\begin{math}
  {\cal M}_F=\{u\in C^1:u'(0)=F(u)\}\subset C^1\mbox{.}
\end{math}
Moreover, for times $t$ before blow-up the map
\begin{math} 
  {\cal M}_F\ni u_0 \mapsto u_t\in{\cal M}_F
\end{math}
is continuously differentiable. Thus, sd-DDEs generate a $C^1$
semiflow (time-$t$ maps) in suitable open subsets of ${\cal M}_F$ (for
example, in a sufficiently small neighborhood of equilibria or
periodic orbits). Hence, Walther's result immediately implies that the
principle of linearized stability applies with respect to
perturbations in ${\cal M}_F$, in particular to equilibria\cite{Stu15}
and periodic orbits. This basic existence result requires only
first-order mild differentiability (a slightly weaker version of them,
since continuity of $F$ in $C^0$ is not needed \cite{W03,
  HKWW06}). Krisztin \cite{K03} proved that the unstable manifold of
equilibria is a $C^\ell$ graph for $\ell$ times mildly differentiable
right-hand sides, using a slightly different (possibly equivalent)
definition of mild differentiability for orders greater than
$1$. Based on Walther's semiflow results,
\chg{Stumpf\cite{Stu11,Stu16b} proved the existence and attractivity
  of $C^1$ local center-unstable and center manifolds near
  equilibria. Alternative proofs are given by
  Krisztin\cite{K06,K06a}.} Furthermore, the assumptions (S1--S3)
imply that periodic boundary-value problems are equivalent to
finite-dimensional smooth systems of algebraic equations for a
sufficiently large number of first Fourier
coefficients\cite{S12}. This equivalence permits us to perform a
classical Lyapunov Schmidt reduction near equilibria $u_*$ for which
the characteristic matrix $\Delta(\lambda)\in \C^{n\times n}$, defined
by $\Delta(\lambda)q=\lambda q-\partial F(u_*)[\theta\mapsto
q\exp(\lambda\theta)]$ has a single pair of roots on the imaginary
axis. Consequently, the classical Hopf bifurcation theorem about a
family of periodic orbits branching off from $u_*$ is
valid\cite{E06,S12}, including formulas determining criticality of the
Hopf bifurcation. More generally, the reduction of periodic boundary
value problems to smooth algebraic equations implies that all objects
computed by DDE-Biftool depend as expected on parameters and the
right-hand side such that they can be computed using standard
numerical discretizations\cite{S12}. This includes branches of
periodic orbits in parameter-dependent systems, the variational
problems for folds, period doublings and torus
bifurcations\cite{S13}. Statements about periodic orbit families
branching off at period doublings and resonant torus bifurcations (in
resonance tongues, first computational demonstrations for DDEs were
for an El-Nin{\~o} model\cite{TCZXB98,KKP15,KKP16,KS14}) follow in a
similar way from a Lyapunov-Schmidt reduction as the Hopf bifurcation
statement.

\section{Normal form computations in DDEs with constant delays --- Review}
\label{sec:nf:const}
\setcounter{paragraph}{0}
Recent work by Kuznetsov, Janssens, Wage and
Bosschaert\cite{J10,W14,B16,BJK1x} has developed and implemented
expressions for the normal form coefficients of local bifurcations in
DDEs with constant delays. For discrete delays, this corresponds to
the case where the delay functions $\tau^j$ in \eqref{fde:tau} are all
constant (e.g., parameters) independent of the state. Their procedure
follows closely the methods originally developed for ODEs\cite{K99}
(and is in principle applicable to other abstract
ODEs\cite{GJKV13}). They assume that the DDE $\dot u(t)=F(u_t)$ has
an equilibrium at $u_*$. For our notation we assume $F(0)=0$, and denote
the first derivative of the right-hand side $F:C^0\to\R^n$ in $0$ by
$A=\partial F(0)\in{\cal L}(C^0;\R^n)$.

\subsection{Linear stability and center manifold}
\label{sec:lin:cmf}
The matrix $\Delta(\lambda)\in \C^{n\times n}$ defined by $\Delta(\lambda)q=\lambda
q-A[\theta\mapsto q\exp(\lambda\theta)]$ for
$q\in\C^n$ is called the \emph{characteristic matrix}. We assume that
the characteristic equation
\begin{align*}
  \det \Delta(\lambda)&=0
\end{align*}
has $n_c$ roots (including multiplicity) on the imaginary axis:
\begin{displaymath}
  \sigma_c=\{\lambda_1,\ldots,\lambda_{n_c}\}=
  \{\lambda\in\C:\det\Delta(\lambda)=0\}\cap \i\R\mbox{.}
\end{displaymath}
 For the type of functionals $F$ that
DDE-Biftool treats, $\Delta(\lambda)$ is given by
\begin{displaymath}
  \Delta(\lambda)=\lambda\id-
  \sum_{j=1}^m\partial_jf(0,\ldots,0)\e^{-\lambda\tau^j}\mbox{,}
\end{displaymath}
where for constant delays the $\tau^j$ are parameters, while for
state-dependent delays, the $\tau^j$ are evaluated at the equilibrium
$0$. The corresponding eigenvectors are in $C^\infty$, and have the
form $\theta\mapsto q\exp(\lambda\theta)$. The generalized
eigenvectors (also in $C^\infty$ if present) have the form
$\theta\mapsto \sum_{j=0}^{j_{\max}}q^j\theta^j\exp(\lambda\theta)$,
where $j_{\max}+1$ is the length of the Jordan chain and
$q^0,\ldots,q^{j_{\max}}$ are in $\C^n$. Let
$B=\{b_1,\ldots,b_{n_c}\}$ be a basis of real functions of the linear
center subspace $U_c=\lspan B$ of $\dot u=Au_t$ in $C^0$, and let
$B^\dagger:C^0\to \R^{n_c}$ be such that $B^\dagger B=\id$ in
$\R^{n_c}$ and $BB^\dagger$ is a spectral projection onto $\lspan B$
(see \eqref{eq:pc}--\eqref{eq:resolvent} in the Appendix for a
concrete expression based on the resolvent formalism).

\paragraph*{Center manifold for constant delays}
For DDEs with constant discrete delays ($\tau^j=\mathrm{const}$ in
\eqref{dom:tau}) the time-$t$ map $C\ni u_0\mapsto u_t\in C$ is as
smooth\cite{DGLW95, HL93} as the right-hand side $f:\R^{n\times
  m}\mapsto \R^n$ in \eqref{dde:sys_rhs}. The reason is that, for
those $f$, the right-hand side as a map $F:C^0\to\R^n$ is
smooth. Hence, in a ball $B_r(0)$ around $0$ with sufficiently small
radius $r$ a smooth center manifold of dimension $n_c$,
$h:B_r(0)\subset\R^{n_c}\to C^0$ exists.

\chg{More precisely, let us assume that the right-hand side coefficient
function $f$ in \eqref{dde:sys_rhs} is at least $\ell$ times continuously
differentiable. Then we can find a radius $r>0$ such that the
invariant graph $h:B_r(0)\subset\R^{n_c}\to C^\ell$ is $\ell$ times
differentiable\cite{DGLW95,HL93}. We write the graph as $h(\theta;u_c)$, putting the
argument of the function $h(u_c)$ in $C^\ell$ first. For any initial
condition $u_0(\theta)=h(\theta;u_c^0)$ ($u_c^0\in B_r(0)$) on the
graph, $u_t(\theta)$ equals $h(\theta;u_c(t))$, where
\begin{equation}\label{eq:cmf:ode}
  \dot u_c(t)=B^\dagger \partial_1h(\cdot;u_c(t))\mbox{,}
\end{equation}
and $u_c(0)=u_c^0$, as long as $|u_c(t)|\leq r$.}

\subsection{Normal form computation}
\label{sec:nf:dde}
Assuming that the right-hand side $F$ and the center manifold $h$ are
smooth up to a desired order $\ell$ (as is the case for constant
delays), it is known that the flow on the local center manifold can be
brought into a normal form up to order $\ell$, such that the flow on
the center manifold $\dot u_c=B^\dagger\partial_1 h(\cdot;u_c)$ has a given expansion
\begin{equation}
  \label{eq:cmf:odenf}
  \dot u_c=A_c^1 u_c+\sum_{j=2}^\ell \frac{1}{j!}A_c^j[\alpha_j]u_c^j+o(|u_c|^\ell)\mbox{.}
\end{equation}
Equation~\eqref{eq:cmf:odenf} is an ODE for $u_c\in\R^{n_c}$. \chg{All
  derivatives up to order $\ell$ of the remainder $o(|u_c|^\ell)$ are
  smaller than the corresponding derivatives of the lower-order terms
  for all small $|u_c|$.} All of the $j$-linear coefficients $A_c^j$
depend only on the type of equilibrium (which local bifurcation?),
except for the still-to-be-determined normal form parameters
$\alpha_j$ at each order $j>1$. The linear coefficients $A_c^1$ are
uniquely determined by $B$ and $B^\dagger$: $A_c^1=B^\dagger B'$,
where $B'$ is the derivative of $B$ with respect to the space variable
$\theta$. \chg{There exists a $C^\ell$-smooth coordinate change in
  $\R^{n_c}$ that transforms the ODE~\eqref{eq:cmf:ode}, describing
  the semiflow of the DDE restricted to its local center manifold $h$,
  into Equation~\eqref{eq:cmf:odenf} (this is called smooth local
  equivalence).}

Normal form computations are concerned with the computations of these
unknown coefficients $\alpha_j$ and, if desired, the expansion
coefficients $h_j(\theta)=\partial^j_2h(\theta;0)$ of the center
manifold. Inputs are the expansion coefficients $F_j=\partial^j F(0)$
(also $j$-linear forms) of the right-hand side of the DDE, and the
general parametric normal form expansion coefficients $A_c^j[\cdot]$,
which depend on the type of the bifurcation investigated (e.g., Hopf
bifurcation and degenerate Hopf bifurcation in the example in
Section~\ref{sec:pos}). The procedure for computing the coefficients
$\alpha_j$, as outlined for ODEs by Kuznetsov\cite{K99}, and adapted
to DDEs recently\cite{J10,W14,B16,BJK1x}, is summarized in
Section~\ref{app:locbif} in the appendix.

The invariance of $h$ gives at each order a linear system of equations
for the expansion coefficients $h_j(0)$ of the center manifold at
$\theta=0$. The system depends also linearly on $\alpha_j$ (if at
order $j$ a normal form coefficient is present). The coefficients of
the linear system for $h_j(0)$ and $\alpha_j$ depend only on $A$ (same
as $F_1$), the linear part of $F$. At each order $j$, the coefficient
$\alpha_j$ is determined by the Fredholm alternative as the unique
value for which the linear system is solvable for $h_j(0)$.

\subsection{General example --- Hopf bifurcation}
\label{sec:hopf}
A typical result of the procedure is the normal form coefficient $L_1$
(which would be the real part of $\alpha_3$, divided by $\omega$) for
the Hopf bifurcation\cite{W14}, as implemented in
DDE-Biftool\cite{ddebiftoolmanual,W14,B16,BJK1x}. Suppose the
linearized DDE $\dot u=\partial F(0)u_t=Au_t$ has a purely imaginary
eigenvalue pair $\pm\i\omega$, with the eigenvector
$q=q_0\e^{\i\omega\theta}$ and its complex conjugate $\bar q=\bar
q_0\e^{-\i\omega}$. That is,
\begin{displaymath}
  \Delta(\i\omega)q_0=\i\omega q_0-A[\e^{\i\omega\theta}q_0]=0\mbox{,}
\end{displaymath}
\chg{and $\pm\i\omega$ are the only roots of $\det \Delta(\cdot)$ on the
imaginary axis}. For notational convenience one chooses as basis
$B=h_1$ of the center subspace of $C^0$ the vectors $\{q,\bar q\}$,
thus using complex notation instead of, for example, $\{\re q,\im
q\}$. The projection $B^\dagger$ is given by the normalized adjoint
eigenvector $p$ for $\i\omega$ and its complex conjugate $\bar p$. The
general expression for adoint eigenvectors is given by Diekmann
\emph{et al}\cite{DGLW95}. For the particular case, where the linear
functional $A$ has the form
\begin{displaymath}
  Au=\sum_{j=1}^mA_ju(-\tau^j)
\end{displaymath}
(as arising in problems treatable with DDE-Biftool) and the critical
spectrum consists of simple eigenvalues $\pm\i\omega$, the projection
is of the form
\begin{align*}
  B^\dagger_1 u&=p_0u(0)+
  \sum_{j=1}^m\int_0^{\tau^j}\e^{\i\omega s}p_0A_ju(s-\tau^j)\d s\mbox{,}\\
  B^\dagger_2u&=\bar B^\dagger_1 u\mbox{.}
\end{align*}
The $C^{1\times n}$ vector $p_0$ is given by $p_0\Delta(\i\omega)=0$
and (after normalization) $p_0\Delta'(\i\omega)q_0=1$.  At order $2$
the linear system for the coefficients of the center manifold is
regular (thus, $\alpha_2$ is empty). Solving it yields
\begin{align*}
  h_2^{11}(\theta)&=2\Delta(0)^{-1}F_2\,q\bar q\mbox{,}&
  h_2^{20}(\theta)&=\Delta(2\i\omega)^{-1}F_2\,qq\,\e^{2\i\omega\theta}
\end{align*}
(the remaining coefficient is $h_2^{02}=\bar h_2^{20}$). At order $3$,
there is a single complex coefficient ($\alpha_3\in\C$ of which the
real part is the coefficient $\omega L_1$) such that:
\begin{equation}\label{eq:hopf:l1}
  L_1=\frac{1}{2\omega}\re\left(p_0\left[F_3\,qq\bar q+F_2\,\bar
    qh_2^{20}+F_2\,qh_2^{11}\right]\right)\mbox{.}
\end{equation}
If the coefficient $L_1$ is non-zero the Hopf bifurcation is
non-degenerate (subcritical if $L_1>0$, supercritical if $L_1<0$).

\section{Extension to DDEs with state-dependent delays}
\label{sec:nf:sd}
Several observations about the normal form reduction imply that at
least the computational procedure can be extended to DDEs with
state-dependent delays (sd-DDEs). 

The procedure described in section~\ref{sec:nf:dde} requires the
expansion coefficients $F_j$ of the nonlinearity $F$ up to the desired
order (often at least $3$). However, we observe that the derivatives
are applied only to deviations that are expansion coefficients of the
center manifold, $(\theta,u_c)\mapsto h_j(\theta)u_c^j$, where
$\theta$ is the history variable and $u_c$ is the deviation along the
center manifold.  At each order $j$, the unknown coefficient
$h_j(\theta)$ is a solution of the linear ODEs \eqref{hcoeff:ode} (see
Appendix) with constant coefficients and an inhomogeneity that is a
linear combination of $h_k(\theta)$ from lower orders ($k<j$). The
basis of the linear center subspace (called $B$ in the previous
section and equal to $h_1$) consists of functions of the form of a
finite sum
\begin{equation}\label{eq:expform}
  \theta\mapsto\sum_{i=1}^{n_{\max}}q_i\theta^{\kappa_i}\re\e^{\lambda_i\theta}
\end{equation}
of some length $n_{\max}$ with $n_{\max}$ non-negative integer powers
$\kappa_i$ of $\theta$ (possibly, some $\kappa_i=0$), and complex
exponents $\lambda_i$. Therefore the ODE \eqref{hcoeff:ode} defining
the coefficients $h_j(\theta)$ implies that all center manifold
expansion coefficients have the form \eqref{eq:expform}. Hence, they
are smooth in $\theta$ such that the functional $F$ can be
differentiated in the equilibrium in the direction of $\sum_{j=1}^\ell
h_j(\theta)u_c^j$ for all $\ell$ and all $u_c\in\R^{n_c}$.

The derivative of expressions of the form \eqref{eq:expform} is known
analytically such that a user routine computing the directional
derivative
\begin{displaymath}
  \frac{\partial^\ell}{\partial \delta^\ell} F\left.\left(\delta 
    \sum_{j=1}^\ell h_j(\theta)u_c^j\right)\right\vert_{\delta=0}
\end{displaymath}
can rely on all derivatives of the argument of $F$ with
respect to $\theta$. Similarly, finite-difference approximations of
the derivative with respect to $\delta$ are known to converge. Both
approaches are experimentally supported in the current development
version of DDE-Biftool\cite{ddebiftoolmanual}. Section~\ref{sec:pos}
will illustrate their use for a position control problem.

\subsection{Illustration for Hopf bifurcation in sd-DDE~\eqref{dde:nested}}
\label{sec:hopf:illu}
For the example $\dot x(t)=p-x(t+x(t))$ the characteristic matrix
$\Delta(\lambda)$ of the linearization in the equilibrium
$x_*=p$ has the form $\Delta(\lambda)=\lambda-\e^{\lambda p}$,
which has a Hopf bifurcation with critical eigenvalue $\i\omega=\i$ at
$p=-\pi/2$. Thus, the right eigenvector is $q(\theta)=\e^{\i\theta}$,
and the left eigenvector $p$ will be scaled such that
$p(0)\Delta'(\i)q(0)=1$. Thus,
$p_0=1/(1+\i\pi/2)\approx\mathtt{0.2884 - 0.4530i}$. The second and
third directional derivatives of $F(u)=p-u(u(0))$ in $0$ along a fixed
direction $v$ are
\begin{align*}
  F_2vv&=-2v(0)v'(-\pi/2)\mbox{,}&
  F_3vvv&=-3v(0)^2v''(-\pi/2)\mbox{.}
\end{align*}
The mixed derivatives $F_2q\bar q$ and $F_3qq\bar q$ can be
constructed from directional derivatives using the polarization
identity (DDE-Biftool's implementation uses this approach).  Following
the procedure for the general Hopf normal form in
Section~\ref{sec:hopf} we compute
$h_2^{20}(\theta)=(0.4+0.8\i)\e^{2\i \theta}$ and
$h_2^{11}(\theta)=-4$ (constant), resulting in a Lyapunov coefficient
\begin{displaymath}
  L_1=\frac{1}{2}\re\left(\frac{2-\i}{1+\i\pi/2}\right)\approx\mathtt{0.0619},
\end{displaymath}
which indicates
that the Hopf bifurcation is subcritical (dangerous) for this example.

\subsection{Smoothness of coefficients}
\label{sec:sd:smooth}

A combination of previous results provides an immediate partial
justification for the normal forms computed with the procedure given
by Kuznetsov \emph{et al}\cite{{J10,W14,B16,BJK1x}} and summarized in
Section~\ref{sec:nf:const}. First of all, trajectories of sd-DDEs
become more regular over time. This effect is well known for DDEs with
constant delays, but also holds for sd-DDEs. \chg{The general proof
requires the precise definition of order-$\ell$ mild
differentiability. We formulate the the statement here for DDEs with
discrete state-dependent delays. }
\begin{proposition}[Smoothness for large times]\label{thm:longtime:smooth}
  \chg{Assume that $F$ is a functional with $C^\ell$ coefficients and $m$
  discrete state-dependent delays \textup{(}of the form
  \eqref{fde:F}--\eqref{fde:tau}\textup{)} less than $\tmax$. } Let
  $u(t)$ with $t\in[-\tmax,T]$ be a solution of $\dot u(t)=F(u_t)$ with
  $u_0\in C^1$ and $u_0'(0)=F(u_0)$.  Then $u_t\in C^\ell$ if $t\geq
  \ell\tmax$. The $\ell$th derivative $u^{(\ell)}$ satisfies a
  (differential) equation of the form
  \begin{equation}
    \label{eq:dul}
    u^{(\ell)}(t)=F^\ell(u_t)\mbox{,}
  \end{equation}
  where $F^\ell$ \chg{has $C^0$ coefficients and $m_\ell=(m+1)^{\ell-1}m$
  discrete delays less than $\ell\tmax$.}
\end{proposition}
\paragraph*{Proof} We show this statement (inductively). For $\ell=1$
the statement follows from the differential equation with $F^1=F$
($f^1=f$ and $m=m_1$). Assume that we have for $t\geq\ell\tmax$
\begin{equation}\label{eq:dul:iass}
  u^{(\ell)}(t)=f^\ell(u^1,\ldots,u^{m_\ell})\mbox{,}
\end{equation}
where $u^j=u(t-\tau^j_\ell(u^1,\ldots,u^{j-1}))$ and all
$\tau_\ell^j\leq \ell\tmax$ (for $\ell=1$, $\tau_1^j=\tau^j$ for
$j=1,\ldots,m$). Thus, for $t\geq(\ell+1)\tmax$ $u_t(\theta)$ is $C^1$
for all $\tau\in[-\ell\tmax,0]$. Consequently, the right-hand side of
\eqref{eq:dul:iass} is differentiable with respect to time for
$t>(\ell+1)\tmax$ (and, hence, the left-hand side). Its derivative is
\begin{align}\nonumber
  u^{(\ell+1)}(t)&=\partial F^\ell(u_t)\dot u_t\\
  \label{flp1def}
  &=\sum_{j=1}^{m_\ell}\partial_jf^\ell(u^1,\ldots, u^{m_\ell})V^j
  \mbox{\ where}\\
  \nonumber
  u^j&=u_t(-\tau_\ell^j)\mbox{\quad for $j=1,\ldots,m_\ell$,}\\
  \nonumber
  (\partial_k)\tau_\ell^j&=(\partial_k)\tau_\ell^j(u^1,\ldots,u^{j-1})\\
  \label{vjdef}
  V^j&=\dot u(t-\tau_\ell^j)
  \left[1-\sum_{k<j}\partial_k\tau_\ell^jV^k\right]\mbox{.}
\end{align}
\chg{For $j=1$ the above expression \eqref{vjdef} for $V^j$ equals $\dot
u(t-\tau_\ell^1)=\dot u(t)$.} We replace $\dot u(t-\tau_\ell^j)$ in
\eqref{vjdef} with $F^1(u_{t-\tau_\ell^j})$ such that
\begin{align*}
  V^j&=f^1(u^{m_\ell+(j-1)m_1+1}\!\!\!,\ldots,u^{m_\ell+jm_1})
\left[1-\sum_{k<j}\partial_k\tau_\ell^jV^k\right]\mbox{,}
\end{align*}
where for $k=1,\ldots,m_1$
\begin{multline*}
  u^{m_\ell+(j-1)m_1+k}=u\left(t-\tau_\ell^j(u^1,\ldots,u^{j-1})\right.\\
  -\left.\tau_1^k(u^{m_\ell+(j-1)m_1+1},\ldots,u^{m_\ell+(j-1)m_1+k-1})\right)\mbox{.}
\end{multline*}
We see that the right-hand side in \eqref{flp1def} is a functional
$F^{\ell+1}$ of the same form as $F^\ell$, but where $f^{\ell+1}$ has
$m_\ell+m_1m_\ell$ arguments such that we have $m_\ell+m_1m_\ell$
delays. Those delays are $\tau_\ell^1$,\ldots, $\tau_\ell^{m_l}$ and
for $j=m_\ell+(i-1)m_\ell+k$ ($i=1\ldots,m_\ell$, $k=1,\ldots,m_1$)
\begin{multline*}
  \tau_\ell^{i,k}=\tau_\ell^j(u^1,\ldots,u^{i-1})\\
  +\tau_1^k(u^{m_\ell+(i-1)m_1+1},\ldots,u^{m_\ell+(i-1)m_1+k-1})\mbox{,}
\end{multline*}
which are all less than $(\ell+1)\tmax$. Hence, $u^{(\ell+1)}$ exists
for $t>(\ell+1)\tmax$ and satisfies
$u^{(\ell+1)}(t)=F^{\ell+1}(u_t)$. (End of proof of
Proposition~\ref{thm:longtime:smooth})

\chg{Since $F^\ell(0)=0$, and the coefficients $f^j$ and $\tau^k_j$
  are still at least $C^1$ for all $j\leq\ell$ (we have differentated
  only $\ell-1$ times), we have for all $u_0\in C^1$ sufficiently
  close to $0$ that
\begin{equation}
  \label{eq:lipbd}
  \|u^{(j)}_t\|_0\leq C_j(t)\|u_0\|_0
\end{equation}
for $t\geq\ell\tmax$ and all $j\leq\ell$ and some constant $C(t)>0$.}

A local center-unstable manifold $h$ is exists and is continuously
differentiable for functionals $F$ with $C^1$ coefficients and
discrete state-dependent delays, according to Stumpf
\cite{Stu11}. \chg{Consequently, if $F(0)=0$ and the critical spectrum
$\sigma_c$ of $\dot u=\partial F(0)u_t$ is not empty, a continuously
differentiable local center manifold $h$ exists, too (applying the
standard local center manifold theorem to the ODE with $C^1$-smooth
coefficients that one obtains by restricting the sd-DDE onto its local
center-unstable manifold, see also Stumpf's or Krisztin's
arguments\cite{K06,K06a,Stu16b}). A simple backwards extension and
Proposition~\ref{thm:longtime:smooth} permit us to conclude that all
elements of the local center manifold $h$ are in $C^\ell$:
\begin{lemma}[Smoothness on center manifold]\label{thm:hsmooth}
  Assume that $F$ is a functional with $C^\ell$ coefficients and
  discrete state-dependent delays \textup{(}of the form
  \eqref{fde:F}--\eqref{fde:tau}\textup{)}, with $F(0)=0$, a center
  subspace $\lspan B$ of $\Delta(\lambda)=\lambda\id-\partial
  F(0)[\theta\mapsto \exp(\lambda\theta)]$ of dimension $n_c$ and a
  continuously differentiable local center manifold $h:B_r(0)\subset
  \R^{n_c}\to C^1$, defined in a ball $B_r(0)$ of radius $r>0$ in
  $\R^{n_c}$, for $\dot u(t)=F(u_t)$.

  Then there exists a constant $C>0$ and a radius $r_\ell>0$ such
  $h(\cdot;u_c)\in C^\ell$ and $\|h(\cdot;u_c^0)\|_\ell\leq C
  \|h(\cdot;u_c^0)\|_0$ for all $u_c\in B_{r_\ell}(0)$.
\end{lemma}
\emph{Proof} Let $L\geq0$ be the Lipschitz constant for the right-hand
side of the ODE on the center manifold $\dot u_c=B^\dagger\partial_1
h(\cdot;u_c)$ on $B_r(0)$ (if necesssary, choose $r$ sufficiently
small such that $L$ exists). Thus, for all $u_c^0\in B_{r_\ell}(0)$
with $r_\ell<r\exp(-\ell\tmax L)$ the solution of $\dot
u_c=B^\dagger\partial_1h(\cdot;u_c)$ starting from $u_c(0)=u_c^0$ does
not leave $B_r(0)$ for times $t$ with $|t|\leq\ell\tmax$. Thus, the
flow map $U_c:[-\ell\tmax,\ell\tmax]\times B_{r_\ell}(0)\ni
(t,u_c^0)\mapsto u_c(t)\in B_r(0)$ is well defined. However, this
implies that, for every $u_c^0\in B_{r_\ell}(0)$, $h(\cdot;u_c^0)$ is
the solution of the DDE $\dot u=F(u_t)$ starting from
$h(\theta;U_c(-\ell\tmax;u_c^0))$. Consequently, by
Proposition~\ref{thm:longtime:smooth}, $h(\cdot;u_c^0)$ is in
$C^\ell$. The relation between the $\|_\cdot\|_\ell$-norm and the
$\|\cdot\|_0$-norm follows then from estimate \eqref{eq:lipbd} and the
Lipschitz constant for $U_c(-\ell\tmax;\cdot)$. (End of proof of
Proposition~\ref{thm:hsmooth})}

Consequently, we can expand at least $F$ in the expression
$F(h(u_c))$, which is present in the normal form expansion. Humphries
\emph{et al}\cite{HCK16} used this fact to demonstrate for their
example how one can expand a sd-DDE near an equilibrium up to order
$\ell$ such that all terms of order $j\leq \ell$ are $j$-linear (and
have, thus, constant delays). The remainder term is of order
$o(\|u_t\|_0^\ell)$ and has state-dependent delays. One incurs delays
of length up to $\ell\tmax$ such that we have the following statement,
generalizing the approach of Humphries \emph{et al}:
\begin{lemma}[Expansion with longer delays]\label{thm:longdelay}
  Let $F$ be a functional with $C^\ell$ coefficients and $m$ discrete
  state-dependent delays $\tau^1$,\ldots,$\tau^m$ \textup{(}of the
  form \eqref{fde:F}--\eqref{fde:tau}\textup{)}. Let $u_0\in C^1$ be
  sufficiently small with $u_0'(0)=F(u_0)$. Then the segments $u_t$
  solving $\dot u(t)=F(u_t)$ satisfy after time $\ell\tmax$ a sd-DDE
  of the form
  \begin{equation}\label{eq:longdelay}
    \dot u(t)=\sum_{j=1}^\ell F_j(u_t)^j+\chg{o(\|u_t\|_\ell^\ell)}\mbox{.}  
  \end{equation}
  The $j$-linear functionals $F_j$ and the remainder map
  $C([-\ell\tmax,0];\R^n)$ into $\R^n$. The expansion products
  $(u_t)^j$ have delays that are sums $\tau^{k_1}+\ldots+\tau^{k_j}$,
  where $\{k_1,\ldots,k_j\}\subseteq\{1,\ldots,m\}$ and all delays are
  evaluated at $u=0$.
\end{lemma}
\paragraph*{Proof}
Since after time $t\geq\ell\tmax$ the solution $u_t$ is $\ell$ times
continuously differentiable, we can expand the functional $F$ in the
equilibrium $0$ and in the direction of $u_t$ to order $\ell$ using
its classical differentiability when restricted to $C^\ell$:
\begin{equation}\label{eq:expansion:uprime}
  \dot u(t)=\sum_{j=1}^\ell\partial^jF(0)
  [u_t,u_t',\ldots,u_t^{(j-1)}]^j+\chg{o(\|u_t\|_{\ell}^\ell)}\mbox{.}
\end{equation}
In expansion~\eqref{eq:expansion:uprime} the $j$-form $\partial^jF(0)$
is continuous only on functions in $C^{j-1}$. To keep track of this
dependence on the derivatives of $u_t$, we include the derivatives
explicitly into the multi-linear arguments in
\eqref{eq:expansion:uprime}. To get an expansion that depends on
$u_t\in C^0([-\ell\tmax,0];\R^n$ (no derivatives, but longer history),
we recursively replace derivatives $u^{(j)}(t)$ by $F^{j}(u_t)$ (as
obtained in Proposition~\ref{thm:longtime:smooth}), followed by
expansions of $F^j(u_t)$. A functional
$F^j:C([-j\tmax,0];\R^n)\to\R^n$ generates also a map $F^j_{j+k}$
from $C([-(j+k)\tmax,0];\R^n)$ into $C([-k\tmax,0];\R^n$ for any
$k\geq 0$ via $F^j_{j+k}(u_t)(\theta)=F^j(u_{t+\theta})$. The
subscript $j+k$ indicates the length of the time interval that
arguments of $F^j_{j+k}$ should have. Thus, after the first
replacement of $u^{(\nu)}_t$ by $F^\nu_{\nu+1}(u_t)$, we have that for
$t\geq \ell\tmax$, $u$ satisfies
\begin{displaymath}
  \dot u(t)=\sum_{j=1}^\ell\partial^jF(0)
  [u_t,F^1_2(u_t),\ldots,F^{j-1}_j(u_t)]^j+o(\|u_t\|_\ell^\ell)\mbox{.}  
\end{displaymath}
At subsequent expansions terms from lower orders will change
expansions at higher orders. It remains to be shown inductively that
eventually all derivatives disappear except for the remainder, and
that the length of the history segments $u_t$ does never exceed
$\ell\tmax$.

Let us make the inductive assumption that a history segment
$u_t^{(j)}$ of length $k\tmax$ shows up at order $(j+1)k\leq \ell$. In
the first inductive step we have $k=1$, $j\in\{1,\ldots,\ell-1\}$ and
orders at which derivatives of $u_t$ appear from $2$ to $\ell$. When
replacing $u_t^{(j)}$ by $F^j_{(j+1)k}(u_t)$ the history interval
increases to $(j+1)k$. Then $F^j_{(j+1)k}(u_t)$ has to be expanded up
to order $\left\lceil\ell/((j+1)k)\right\rceil$ ($\lceil r\rceil$ is
the lowest integer greater or equal than $r$). In this expansion, we
have $\nu$-linear forms containing derivatives of $u_t$ up to order
$\nu-1$. A derivative of order $i\leq \nu-1$ shows up for orders of
$u_t$ greater or equal than $(i+1)(j+1)k$.

Hence, a term $u^{(j)}_t$ at order $(j+1)k\leq \ell$ creates new $i$th
derivative terms ($i\geq1$) only at order greater or equal than
$(i+1)(j+1)k$ such that the recursion must terminate. (We restrict to
orders less or equal than $\ell$.) Also, the length of the history
interval of the new $i$th derivative term is $(j+1)k\tmax$, which is
less than $\ell\tmax$, since $(j+1)k\leq\ell$ by inductive
assumption.

(End of proof of Lemma~\ref{thm:longdelay})



\chg{We combine the result of Lemma~\ref{thm:hsmooth} with
  Lemma~\ref{thm:longdelay} to sharpen the estimate for solutions of
  the FDE $\dot u(t)=F(u_t)$ starting on the local center manifold:
  $u_0=h(\cdot;u_c^0)$ with $u_c^0\in B_{r_\ell}(0)$. Then the
  remainder term is also of order $o(\|u_t\|_0^\ell)$ (since
  Lemma~\ref{thm:hsmooth} provides an estimates for
  $\|h(\cdot;u_c)\|_\ell$ in terms of $\|h(\cdot;u_c)\|_0$:
\begin{equation}
  \label{eq:hcut}
  \dot u(t)=\sum_{j=1}^\ell F_j(u_t)^j+o(\|u_t\|_0^\ell)\mbox{.}
\end{equation}
Since $u_t=h(\cdot;u_c(t))$, we may also also replace the remainder by
$o(|u_c(t)|^\ell)$. The truncated DDE~\eqref{eq:hcut} (dropping the
remainder term) has only constant delays. Hence, the semiflow and
local center manifold $h_\mathrm{trunc}$ of the truncated
DDE~\eqref{eq:hcut} are smooth, and can, thus, be transformed into
normal form with the procedure described in Section
\ref{sec:nf:dde}. Since this normal form transformation up to order
$\ell$ is independent of terms of order $o(|u_c|^\ell)$ and keeps
these terms at order $o(|u_c|^\ell)$, we have that for $u$ on the
local center manifold $h$ of the non-truncated sd-DDE $\dot
u(t)=F(u_t)$, the center component $u_c=B^\dagger u_t$ satisfies an
ODE equal to the normal form of the truncated DDE~\eqref{eq:hcut}
except for a different remainder (still of order $o(|u_c|^\ell)$). The
result has the form (compare \eqref{eq:cmf:odenf})
\begin{equation}
  \label{eq:cmf:sddde}
  \dot u_c=A_c^1 u_c+\sum_{j=2}^\ell \frac{1}{j!}A_c^j[\alpha_j]u_c^j+o(|u_c|^\ell)\mbox{,}
\end{equation}
where all coefficients $\alpha_j$ are identical to those of the
normal form of the truncated DDE~\eqref{eq:hcut}. However, in contrast to the
constant-delay DDE, only the first derivative of the remainder
$O(|u_c|^\ell)$ is guaranteed to be small for all small $u_c$, but
\emph{not} the higher-order derivatives.} This was also demonstrated
numerically by Humphries \emph{et al}\cite{HCK16} for their
example. Any phenomenon predicted by the normal form that persists
under perturbations of size $o(|u_c|^\ell)$ will also be present in
the sd-DDE. This includes all periodic orbits and their changes of
stability.

\paragraph{Normally hyperbolic invariant manifolds}
For some bifurcations the normal form of the truncated system may
predict the presence of, for example, invariant tori that branch off
along torus bifurcation curves, away from strong resonances ($1:1$ to
$1:4$, see\cite{K04}). Their degree of normal hyperbolicity is
proportional to their distance from the torus bifurcation in the
truncated system. Our perturbation (the remainder term
$o(|u_c|^\ell)$) is $C^1$ small in a ball around $0$, but not
guaranteed to be $C^j$ small compared to lower order terms (with
$j>1$), except in $0$, because the local center manifold has not been
proven to be smooth. Hence, close to the torus bifurcation the
invariant tori may be altered by the remainder term. However, the
region around the torus bifurcation where the invariant tori are not
sufficiently normally hyperbolic shrinks as we approach the
neighborhood of $0$ if the remainder term decreases faster than the
normal hyperbolicity. This is the case if one chooses $\ell$
sufficiently large. For example, Humphries \emph{et al}\cite{HCK16}
indeed reported invariant tori branching off from the torus
bifurcation near the Hopf-Hopf interactions as predicted by the normal
form. In their paper the authors compared for their example the
results from the direct normal form expansion for the sd-DDE as
explained in general in Section~\ref{sec:nf:sd} to the results from
the constant-delay DDE as constructed via Lemma~\ref{thm:longdelay}
and found agreement up to numerical round-off errors.

\section{Illustration - position control}
\label{sec:pos}
\begin{figure}[ht]
  \centering
  \includegraphics[scale=0.8]{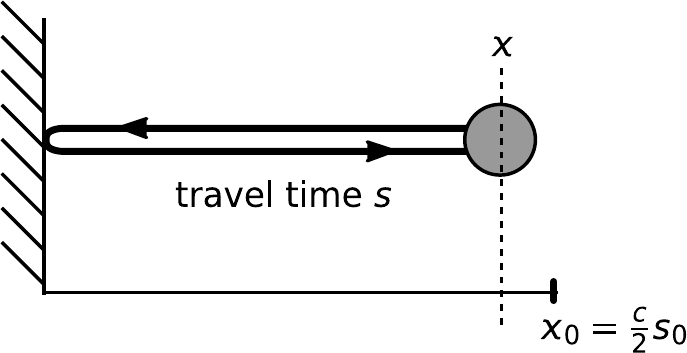}
  \caption{Sketch for position control problem: $x$ is the current
    position of the mover; $x_0$ is the reference position, $c$ is the
    traveling speed of the signal; $s_0$ is the traveling time of the
    signal from obstacle to reference point $x_0$.}
  \label{fig:poscontrol}
\end{figure}

\begin{figure}[ht]
  \centering
  \includegraphics[width=\columnwidth]{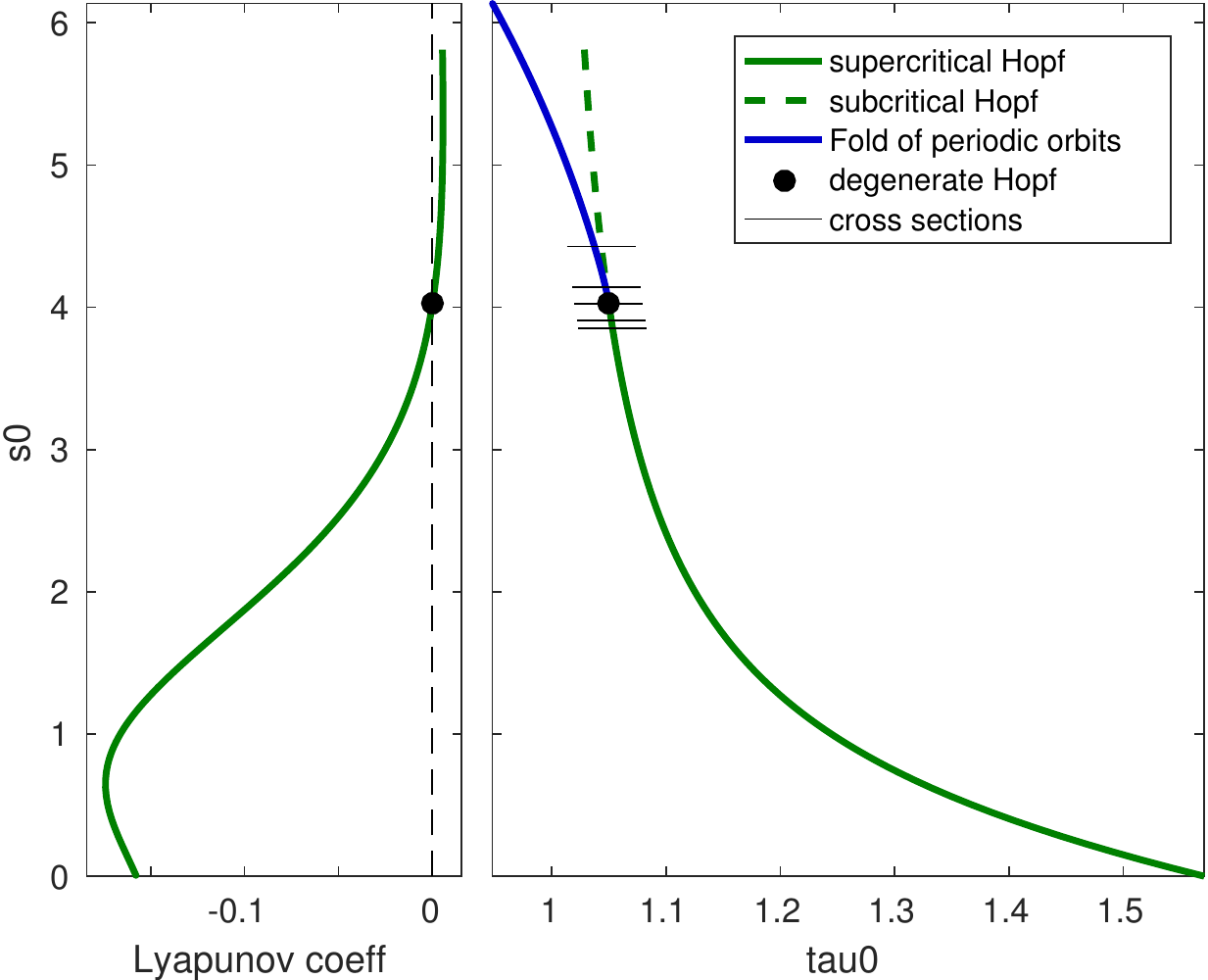}
  \caption{Bifurcation diagram of equilibria and emerging periodic
    orbits in the $(\tau_0,s_0)$-plane, showing the Hopf bifurcation
    and a fold (saddle-node) of periodic orbits. Other parameters: $k=1$,
    $c=2$, $\gamma=1$. Computed with
    DDE-Biftool~\cite{ELR02,ELS01,ddebiftoolmanual} and its normal
    form extension~\cite{W14,B16,BJK1x}.}
  \label{fig:bif2d}
\end{figure}

\begin{figure}[ht]
  \centering
  \includegraphics[width=\columnwidth]{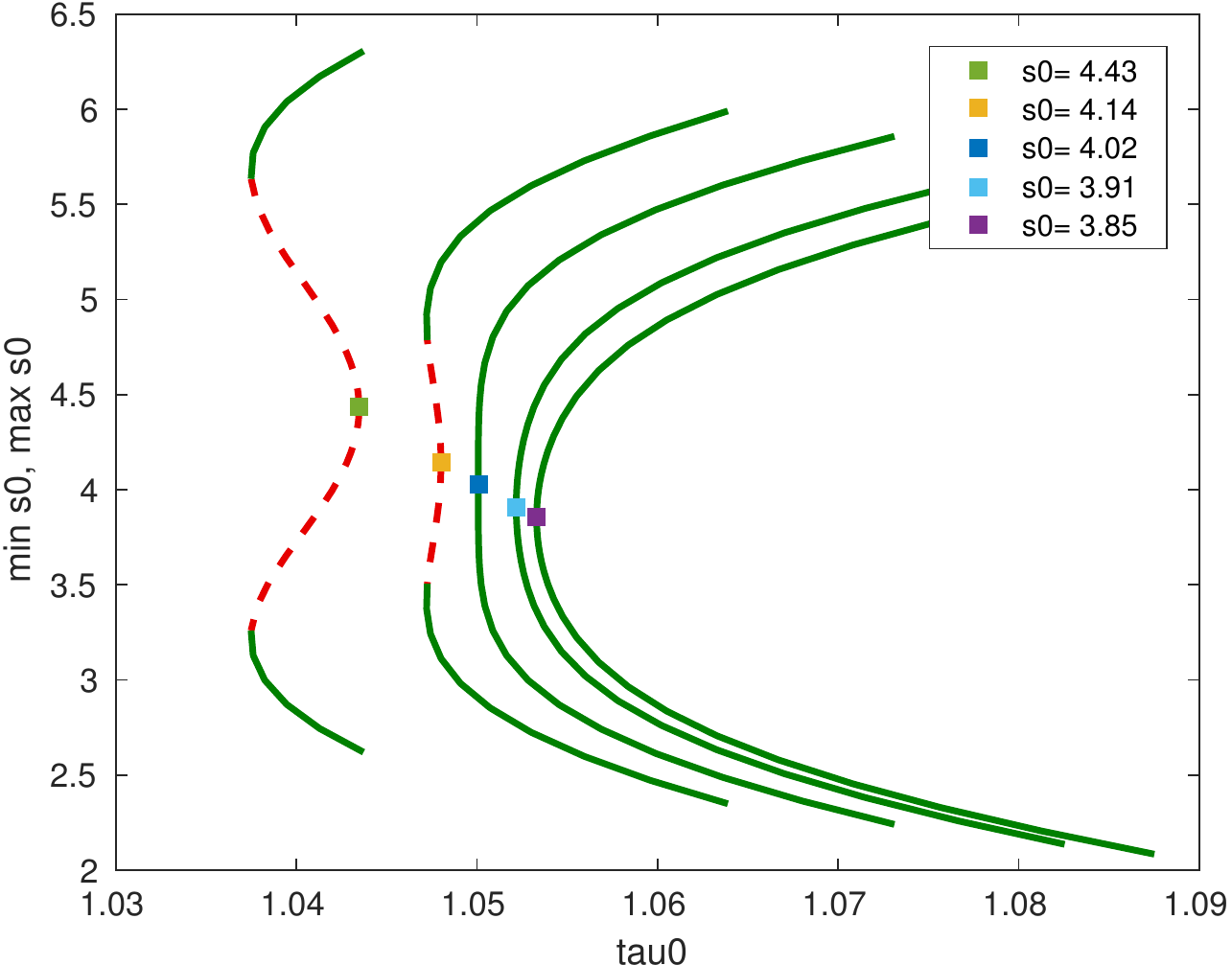}
  \caption{One-parameter families of periodic orbits along the cross
    sections of Figure~\ref{fig:bif2d}: the figure shows maxima and
    minima of the periodic orbits for each parameter value for which
    they have been computed. Dashed curves are unstable periodic
    orbits, solid curves are stable periodic orbits. Other parameters:
    $k=1$, $c=2$, $\gamma=1$. The equilibria undergoing Hopf
    bifurcations are indicated as colored squares. Computed with
    DDE-Biftool~\cite{ELR02,ELS01,ddebiftoolmanual}.}
  \label{fig:bif1d}
\end{figure}

A good example suitable for illustration of simple nonlinear behaviour
introduced by state-dependence of the delay is the position control
problem discussed by Walther \cite{W02} (see also review
\cite{HKWW06}). A mover aims to control its position $x$ relative to
an obstacle using linear position feedback (see Figure~\ref{fig:poscontrol}).
We assume that the controlled motion is free of inertia such that (in
non-dimensionalized quantities)
\begin{equation}\label{eq:control}
  \dot x=k[x_0-x_\mathrm{est}(t-\tau_0)]\mbox{.}
\end{equation}
In \eqref{eq:control} $k$ is the linear control gain, $x_0$ is the
reference position that the mover aims to maintain, $x_\mathrm{est}$
is the mover's estimate of the current position, and $\tau_0$ is a
processing or reaction delay in the control loop. Even if the estimate
$x_\mathrm{est}(t)$ is perfect (equal to $x(t)$), the equilibrium
$x_0$ of the controlled system \eqref{eq:control} will be linearly
unstable if $k\tau_0>\pi/2$. If the mover estimates the current
position by sending out a signal and measuring the traveling time for
the reflected signal then an additional state-dependent delay is
introduced. Let $s(t)$ be the time that the reflected signal, arriving
at the mover time $t$, needed since leaving the mover, and let $c$ be
the signal traveling speed. Then
\begin{equation}
  \label{eq:travel}
  cs(t)=x(t-s(t))+x(t)\mbox{.}
\end{equation}
The mover estimates its current position via
\begin{equation}
  \label{eq:stoxest}
  x_\mathrm{est}=\frac{c}{2}s(t)\mbox{.}
\end{equation}
Let us introduce the reference travel time $s_0=\frac{c}{2}x_0$
corresponding to the reference position $x_0$. The full equation of
motion is
\begin{align}
  \label{eq:x:sd-dde}
  \dot x(t)=&\frac{kc}{2}[s_0-s(t-\tau_0)]\mbox{,}\\
  \dot s(t)=&\frac{2s_0-s(t-\tau_0-s(t))-s(t-\tau_0)}{
    \frac{2}{k}+s_0-s(t-\tau_0-s(t))}\nonumber\\
  &-\gamma\,\frac{cs(t)-x(t)-x(t-s(t))}{c+\frac{kc}{2}[s_0-s(t-\tau_0-s(t))]}\mbox{.}
  \label{eq:s:sd-dde}
\end{align}
The differential equation for $s$ follows from \eqref{eq:control} and
\eqref{eq:travel} via Baumgarte regularization: we rewrite
\eqref{eq:travel} in the form $g(t)=0$ (where
$g(t)=cs(t)-x(t-s(t))-x(t)$), and then replace it by the condition
$\frac{\d}{\d t}g(t)=-\gamma g(t)$, re-arranged for $\dot s(t)$. Every
orbit of \eqref{eq:x:sd-dde}--\eqref{eq:s:sd-dde} that is periodic or
lies on a local center manifold with internal contraction rate less
than $\gamma$ satisfies also the algebraic constraint
\eqref{eq:travel}. When writing system
\eqref{eq:x:sd-dde}--\eqref{eq:s:sd-dde} in the general form $\dot
u=F(u_t)$, the right-hand side of
\eqref{eq:x:sd-dde}--\eqref{eq:s:sd-dde} corresponds to a functional $F$ with the form
($u=(u_1,u_2)^T=(x,s)^T$)
\begin{align*}
  F(u)=
  \begin{bmatrix}
    \frac{kc}{2}[s_0-u_2(-\tau_0)]\\[2ex]
    \begin{matrix}
      \cfrac{2s_0-u_2(-\tau_0-u_2(0))-u_2(-\tau_0)}{
        \frac{2}{k}+s_0-u_2(-\tau_0-u_2(0))}\qquad\\
      \qquad-\gamma\,\cfrac{cu_2(0)-u_1(0)-u_1(-u_2(0))}{c+\frac{kc}{2}[s_0-u_2(-\tau_0-u_2(0))]}
    \end{matrix}
  \end{bmatrix}\mbox{.}
\end{align*}
Equilibria and periodic orbits computed in this illustration had their
$s(t)$ component in the range $[s_{\min},s_{\max}]$ with
$s_{\min}\geq0$ and $s_{\max}<10$ in the parameter ranges used for
figures \ref{fig:bif2d} and \ref{fig:bif1d}. Hence, we may set
$\tmax=10$ and treat $F$ as a functional from $C([-\tmax,0];\R^2)$ to
$\R^2$.

For our demonstration we fix $k=1$, $c=2$ and $\gamma=1$ in
non-dimensionalized quantities. We vary $\tau_0$ and $s_0$ in a
two-parameter bifurcation study. The system has one constant delay
$\tau_0$ and two state-dependent delays. In the notation of
DDE-Biftool the function $f:\R^{2\times 4}\times\R^2\to\R^2$ has
the time-dependent arguments $u(t-\tau_j)=[x(t-\tau_j),s(t-\tau_j)]^T$ for
$j=1,2,3,4$, and the parameters $(\tau_0,s_0)$, where
\begin{align*}
\tau_1&=0\mbox{,}&
\tau_2&=\tau_0\mbox{,}\\
\tau_3&=u_2(t)=s(t)\mbox{,}&
\tau_4&=\tau_0+u_2(t)=\tau_0+s(t)\mbox{.}
\end{align*}
The system~\eqref{eq:x:sd-dde}--\eqref{eq:s:sd-dde} has a unique
equilibrium at $u_*=(x_*,s_*)=(cs_0/2,s_0)$. As part of the principle of
linearized stability proved by Walther\cite{W03} comes the description
for how to compute stability (which is implemented in DDE-Biftool):
``\emph{freeze}'' the state-dependent delays at the values in the
equilibrium, and then compute the linearization of the corresponding
DDE with constant delays
\cite{hartung2001linearized,cooke1996problem,HKWW06}. For the position
control problem this procedure gives a algebraic relation between the
parameter values at which Hopf bifurcations occur:
\begin{align}
  \label{eq:hopf:formula}
  0&=\frac{2\omega^\pm_\ell}{k}-\sin(\omega_\ell^\pm\tau_0)-
  \sin(\omega_\ell^\pm(\tau_0+s_0))\mbox{,\quad where}\\
  \omega_\ell^\pm&=\frac{\pi(1+2\ell)}{\tau_0+s_0\pm\tau_0}\mbox{.}\nonumber
\end{align}
The Hopf bifurcation that forms the boundary of the stability region
in the $(\tau_0,s_0)$-plane is the curve for $\omega_0^+$, shown in
Figure~\ref{fig:bif2d} (right panel) as a green dashed/solid curve.
As expression~\eqref{eq:hopf:formula} is still implicit, the curve in
Figure~\ref{fig:bif2d} was computed with DDE-Biftool. The standard
Hopf bifurcation theorem can be applied to sd-DDEs\cite{E06,S12} such
as system~\eqref{eq:x:sd-dde}--\eqref{eq:s:sd-dde}. Hence, a family of
periodic orbits branches off from the Hopf bifurcation. Near the
equilibrium the stability of periodic orbits can be predicted using
the expression \eqref{eq:hopf:l1} for $L_1$ as implemented by
Kuznetsov \emph{et al}\cite{W14,B16,BJK1x}. This was
rigorously proven using a Lyapunov-Schmidt reduction for periodic
boundary value problems\cite{S12}. Its value along the Hopf curve is
shown in the left panel of Figure~\ref{fig:bif2d}. The value of $L_1$
crosses zero at $s_0\approx 4.02$, $\tau=1.05$. There the Hopf
bifurcation is degenerate and the second Lyapunov coefficient is
$L_2\approx -1.9\times 10^{-3}$. This implies that the family of
periodic orbits exists to the right and is stable where the Hopf
curve is solid in Figure~\ref{fig:bif2d}. The family of periodic
orbits is unstable and exists to the left, before folding in a fold of
periodic orbits to the right where the Hopf curve is dashed in
Figure~\ref{fig:bif2d}.

\section{Conclusion}
\label{sec:conc}
As this paper shows, expressions for normal form coefficients for
constant-delay DDEs can be generalized to sd-DDEs. The mathematical
justification is only partially complete, but for many phenomena it is
already clear how they persist when the truncation is removed. The
complete justification requires smoothness for the local center
manifold. Krisztin has provisional results\cite{K06} that show how his
proof for smooth unstable manifolds of equilibria\cite{K03} can be
extended to local center manifolds. Ideally, the general result for
persistence of compact normally hyperbolic manifolds should in some
sense be adapted to sd-DDEs in the following form.  Consider a sd-DDE
of the form
\begin{equation}\label{eq:gendde}
  \dot  u(t)=F_\mathrm{c}(u_t)+F_\mathrm{sd}(u_t),
\end{equation}
where $F_\mathrm{c}: C^0\to \R^n$ is smooth and $\dot
u(t)=F_\mathrm{c}(u_t)$ has a compact overflowing invariant normally
hyperbolic (say, stable) manifold ${\cal M}_0$. If we also assume that
$F_\mathrm{sd}$ has a sufficiently small Lipschitz constant with
respect to the space of Lipschitz continuous functions $C^{0,1}$ (and
is mildly differentiable up to order $\ell$), then \eqref{eq:gendde}
should also have a compact overflowing invariant normally stable
manifold ${\cal M}$. The smoothness of ${\cal M}$ should only be
restricted by the spectral gap in the exponential dichotomy on ${\cal
  M}_0$.

\begin{acknowledgments}
  J.S. gratefully acknowledges the financial support of the EPSRC via
  grants EP/N023544/1 and EP/N014391/1. J.S. has also received funding
  from the European Union's Horizon 2020 research and innovation
  programme under Grant Agreement number 643073.
\end{acknowledgments}
\appendix

\section{Details of normal form expansion for local bifurcations of DDEs}
\label{app:locbif}
This appendix gives a few additional details for the computation of
coefficients in the normal form procedure of
Section~\ref{sec:nf:const}. 
\paragraph*{The linear DDE $\dot u=Au_t$}
Recall that the characteristic matrix is denoted by
$\Delta(\lambda)\in \C^{n\times n}$, which has $n_c$ eigenvalues on
the imaginary axis (counting multiplicity). Let
$B=\{b_1,\ldots,b_{n_c}\}$ be a basis of the linear center subspace
$U_c=\lspan B$ of $\dot u(t)=Au_t$. A spectral projection $P_c$ onto
the space $U_c$ is given by residue of the resolvent $R(\lambda)$:
\begin{align}
  P_c:&C^0\to U_c=\lspan B\mbox{,}&
  P_cv&
=\frac{1}{2\pi\i}\oint_{\textstyle\sigma_c} R(\lambda)
\d \lambda\, v\label{eq:pc}
\end{align}
where the curve integral is taken around the critical spectrum
$\sigma_c$. The resolvent $R(\lambda)$, mapping $C^0$ into
$C^1$ is defined as the unique solution $x\in C^1$ of 
\begin{align*}
  \begin{bmatrix}
    v(0)\\ v(\theta)
  \end{bmatrix}=
  \begin{bmatrix}
    \lambda x(0)-Ax\\
    \lambda x(\theta)-x'(\theta)
  \end{bmatrix}\mbox{,}
\end{align*}
which is
\begin{align}
  x(\theta)=\e^{\lambda\theta}x_0+\int_\theta^0\e^{\lambda(\theta-s)}v(s)\d
    s\mbox{, where}\label{eq:resolvent}\\
    x_0=\Delta(\lambda)^{-1}\left[v(0)+A\left[\int_\theta^0\e^{\lambda(\theta-s)}v(s)\d
    s\right]\right]\nonumber
\end{align}
\chg{We define $B^\dagger:C^0\ni x\mapsto x_c\in\R^{n_c}$, where
$x_c\in\R^{n_c}$ is the unique vector of coordinates such that
$Bx_c=P_cx$. Thus, $B^\dagger B$ is the identity in $\R^{n_c}$, and
$BB^\dagger=P_c$.}

\paragraph*{Center manifold expansion}
The semiflow of the DDE, restricted to the center manifold $\{u\in
C^0: u(\theta)=h(\theta;u_c), u_c\in\R^{n_c}\mbox{\ small}\}$,
introduced in Section~\ref{sec:nf:const}, satisfies the ODE in $\R^{n_c}$
\begin{equation}\label{app:cmf:ode}
  \dot u_c=B^\dagger\partial_1 h(\cdot;u_c)\mbox{.}
\end{equation} The invariance of graph of the
manifold
\begin{displaymath}
  \R^{n_c}\supset B_r(0)\ni u_c\mapsto h(\cdot;u_c)\in C^\ell
\end{displaymath}
under the DDE $\dot u=F(u_t)$ implies
\begin{align}\label{eq:inv0}
  \partial_1 h(0;u_c)&=F(h(u_c))\mbox{, and for $\theta\in[-\tau,0]$}\\
  \label{eq:inv:theta}
  \partial_1 h(\theta;u_c)&=\partial_2h(\theta;u_c)\, \dot
  u_c
  \mbox{.}
\end{align}
Let us introduce expansions for $F$ and $h(\theta;\cdot)$ up to order
$\ell$ in the point $u=0$ (for $F$) and $u_c=0$ (for $h(\theta,\cdot)$):
\begin{align*}
  h(\theta;u_c)&=\sum_{j=1}^\ell \frac{1}{j!}h_j(\theta) [u_c]^j+O(|u_c|^{\ell+1})\mbox{,} \\
  F(u)&=\sum_{j=1}^\ell  \frac{1}{j!}F_j[u]^j+O(|u_c|^{\ell+1})\mbox{.}
\end{align*}
The first-order coefficient $F_1$ of $F$ is the linear operator $A$,
the first-order coefficient $h_1(\theta)$ of the manifold graph is
$B(\theta)$. The coefficients $h_j$ for $j>1$ are only determined up to
conjugacy of the flow on the center manifold to order $j$. A different
choice of $h_j$ corresponds to a different, but conjugate, ODE for
$u_c$. For example, requiring $B^\dagger h_j [u_c]^j=0$ for all $j>1$
and all $u_c\in\R^{n_c}$ would determine $h_j$ uniquely in combination
with the invariance \eqref{eq:inv0}--\eqref{eq:inv:theta}. 

\paragraph*{Determining systems for coefficients $h_j(0)$ and $\alpha_j$}
However,
the approach proposed by Kuznetsov\cite{K99} and taken in DDE-Biftool's normal
form extension\cite{J10,W14,B16,BJK1x} is to choose the expansion
coefficients $h_j$ such that the ODE \eqref{app:cmf:ode} on the center
manifold for $u_c$ is already in normal form:
\begin{equation}
  \label{app:cmf:odenf}
  \dot u_c=A_c^1 u_c+\sum_{j=2}^\ell \frac{1}{j!}A_c^j[\alpha_j][u_c]^j+O(|u_c|^\ell)\mbox{.}
\end{equation}
In \eqref{app:cmf:odenf} the matrix $A_c^1=B^\dagger \circ
[\partial/\partial\theta]\circ B=B^\dagger\circ B'\in \R^{n_c\times
  n_c}$ is the projection of the linear DDE on the eigenspace for the
spectrum $\sigma_c$ on the imaginary axis. For higher orders $j>1$ the
coefficients $A_c^j$ are given except for a finite number of
to-be-determined normal form coefficients $\alpha_j$. We use square
brackets to indicate that $A_c^j$ is a given map depending linearly on
$\alpha_j$ and $j$-linearly on $u_c$. The coefficient $\alpha_j$ may
be empty (for example, $\alpha_1$ is always empty). Inserting the
expansions for $h$, $F$ and $\dot u_c$ into the invariance equation
\eqref{eq:inv:theta} gives at order $j$ a $n$-dimensional
inhomogeneous constant-coefficient differential equation for each
coefficient of the symmetric $j$-form $h_j(\theta)$:
\begin{equation}\label{hcoeff:ode}
  \begin{split}
    h_j'(\theta)[u_c]^j=&\ jh_j(\theta)[u_c]^{j-1}[A_c^1u_c]\\
    &\ +B(\theta)A_c^j[\alpha_j][u_c]^j+R_j(\theta)[u_c]^j\mbox{,}
  \end{split}
\end{equation}
where
\begin{displaymath}
  R_j(\theta)[u_c]^j=\frac{1}{j+1}\sum_{k=2}^{j-1}{j+1 \choose k}h_k(\theta)[u_c]^{j-k}[A_c^k[u_c]^k]
\end{displaymath}
is a known function determined by orders lower than $j$ (it is not
present for orders $1$ and $2$. Let us denote the solution $h_j$ of
the affine ordinary differential equation \eqref{hcoeff:ode} by
\begin{displaymath}
  [H_j(\theta)h_j^0+H_{\alpha,j}(\theta)\alpha_j+
  H_{R,j}(\theta)][u_c]^j\mbox{.}
\end{displaymath}
The above expression indicates that the solution is linear in
$h_j^0=h_j(0)$ (its initial value), $\alpha_j$ and $R_j$, and
$j$-linear in $u_c$.  If the basis $B$ consists only of eigenvectors
(eigenvector $b_i$ for eigenvalue $\lambda_i$), then $A_c^1$ is
diagonal, and $H_j(\theta)=\exp(\lambda_i\theta)h_{j,\nu}^0$ for
coefficients $h_{j,\nu}$ of the $j$-form $h_j(\theta)$. In this case
the ${n+j \choose j}$ differential equations for the ${n+j \choose j}$
coefficients $h_{j,\nu}$ of the $j$-form $h_j(\theta)$ decouple. The
initial conditions $h_j^0$ are determined by the invariance at
$\theta=0$, \eqref{eq:inv0}:
\begin{align*}
  h_j'(0)[u_c]^j&=[Ah_j(\cdot)][u_c]^j+R_j^F[u_c]^j\mbox{,\quad where}\\
  R_j^F[u_c]^j&=\sum_{k=2}^j\sum_{\nu\in\ind(k,j)}F_k
  \prod_{\mu=1}^k h_{\nu_\mu}[u_c]^{\nu_\mu}\mbox{.}
\end{align*}
The second sum is taken over multi-indices $\nu\in\ind(k,j)$. The set
$\ind(j,k)$ is the set of $k$-tuples of positive integers summing up
to $j$.  Inserting the differential equation for $h_j$ and its
solution $H_j$ at $\theta=0$ results in an affine equation for
$h_j^0$ and $\alpha_j$ (the \emph{homological equation}):
\begin{align}
  [L_{h,j}h_j^0][u_c]^j=&\ [L_{\alpha,j}\alpha_j][u_c]^j\label{eq:homol}\\
  &\ +[R_j(0)-R_j^F-AH_{R,j}(\cdot)][u_c]^j\nonumber\\
  \intertext{where}  \nonumber
  [L_{h,j}h_j^0][u_c]^j=&\ \left[AH_j(\cdot)h_j^0\right][u_c]^j- jh_j^0[u_c]^{j-1}[A_c^1u_c]\\
  \nonumber
  [L_{\alpha,j}\alpha_j][u_c]^j=&\ B(0)A_c^j[\alpha_j][u_c]^j-
  \left[AH_{\alpha,j}(\cdot)\alpha_j\right][u_c]^j
\end{align}
One can determine $h_j^0$ and $\alpha_j$ for each $j$ by comparing
coefficients of this $j$-form in $u_c$. For orders $j$, for which the
square coefficient matrix $L_{h,j}$ is regular, the normal form
coefficient $\alpha_j$ is not present (since all terms at this order
are non-resonant). If the matrix $L_{h,j}$ is singular with kernel
dimension $d_j$, then the dimension of $\alpha_j$ is $d_j$ and the
dependence of $A_c^j$ on $\alpha_j$ is such that
$[L_{h,j},-L_{\alpha,j}]$ has full rank. Thus, there is a unique
coefficient $\alpha_j$, for which \eqref{eq:homol} is solvable for
$h_0^j$. The solution $h_0^j$ is not unique, but can be made unique,
for example, by forcing it to be orthogonal to the nullspace of
$L_{h,j}^T$; see the references\cite{J10,W14,B16,BJK1x}.

\bibliography{delay}

\begin{thebibliography}{42}%
\makeatletter
\providecommand \@ifxundefined [1]{%
 \@ifx{#1\undefined}
}%
\providecommand \@ifnum [1]{%
 \ifnum #1\expandafter \@firstoftwo
 \else \expandafter \@secondoftwo
 \fi
}%
\providecommand \@ifx [1]{%
 \ifx #1\expandafter \@firstoftwo
 \else \expandafter \@secondoftwo
 \fi
}%
\providecommand \natexlab [1]{#1}%
\providecommand \enquote  [1]{``#1''}%
\providecommand \bibnamefont  [1]{#1}%
\providecommand \bibfnamefont [1]{#1}%
\providecommand \citenamefont [1]{#1}%
\providecommand \href@noop [0]{\@secondoftwo}%
\providecommand \href [0]{\begingroup \@sanitize@url \@href}%
\providecommand \@href[1]{\@@startlink{#1}\@@href}%
\providecommand \@@href[1]{\endgroup#1\@@endlink}%
\providecommand \@sanitize@url [0]{\catcode `\\12\catcode `\$12\catcode
  `\&12\catcode `\#12\catcode `\^12\catcode `\_12\catcode `\%12\relax}%
\providecommand \@@startlink[1]{}%
\providecommand \@@endlink[0]{}%
\providecommand \url  [0]{\begingroup\@sanitize@url \@url }%
\providecommand \@url [1]{\endgroup\@href {#1}{\urlprefix }}%
\providecommand \urlprefix  [0]{URL }%
\providecommand \Eprint [0]{\href }%
\providecommand \doibase [0]{http://dx.doi.org/}%
\providecommand \selectlanguage [0]{\@gobble}%
\providecommand \bibinfo  [0]{\@secondoftwo}%
\providecommand \bibfield  [0]{\@secondoftwo}%
\providecommand \translation [1]{[#1]}%
\providecommand \BibitemOpen [0]{}%
\providecommand \bibitemStop [0]{}%
\providecommand \bibitemNoStop [0]{.\EOS\space}%
\providecommand \EOS [0]{\spacefactor3000\relax}%
\providecommand \BibitemShut  [1]{\csname bibitem#1\endcsname}%
\let\auto@bib@innerbib\@empty
\bibitem [{\citenamefont {Janssens}(2010)}]{J10}%
  \BibitemOpen
  \bibfield  {author} {\bibinfo {author} {\bibfnamefont {S.~G.}\ \bibnamefont
  {Janssens}},\ }\emph {\bibinfo {title} {On a normalization technique for
  codimension two bifurcations of equilibria of delay differential
  equations}},\ \href@noop {} {\bibinfo {type} {Master's thesis}},\ \bibinfo
  {school} {Utrecht University (NL)} (\bibinfo {year} {2010}),\ \bibinfo {note}
  {supervised by Y.A. Kuznetsov and O. Diekmann}\BibitemShut {NoStop}%
\bibitem [{\citenamefont {Wage}(2014)}]{W14}%
  \BibitemOpen
  \bibfield  {author} {\bibinfo {author} {\bibfnamefont {B.}~\bibnamefont
  {Wage}},\ }\emph {\bibinfo {title} {Normal form computations for Delay
  Differential Equations in {DDE-Biftool}}},\ \href@noop {} {\bibinfo {type}
  {Master's thesis}},\ \bibinfo  {school} {Utrecht University (NL)} (\bibinfo
  {year} {2014}),\ \bibinfo {note} {supervised by Y.A. Kuznetsov}\BibitemShut
  {NoStop}%
\bibitem [{\citenamefont {Bosschaert}(2016)}]{B16}%
  \BibitemOpen
  \bibfield  {author} {\bibinfo {author} {\bibfnamefont {M.~M.}\ \bibnamefont
  {Bosschaert}},\ }\emph {\bibinfo {title} {Switching from codimension 2
  bifurcations of equilibria in delay differential equations}},\ \href@noop {}
  {\bibinfo {type} {Master's thesis}},\ \bibinfo  {school} {Utrecht University
  (NL)} (\bibinfo {year} {2016}),\ \bibinfo {note} {supervised by Y.A.
  Kuznetsov}\BibitemShut {NoStop}%
\bibitem [{\citenamefont {Bosschaert}, \citenamefont {Janssens},\ and\
  \citenamefont {Kuznetsov}(2017)}]{BJK1x}%
  \BibitemOpen
  \bibfield  {author} {\bibinfo {author} {\bibfnamefont {M.~M.}\ \bibnamefont
  {Bosschaert}}, \bibinfo {author} {\bibfnamefont {S.~G.}\ \bibnamefont
  {Janssens}}, \ and\ \bibinfo {author} {\bibfnamefont {Y.~A.}\ \bibnamefont
  {Kuznetsov}},\ }\bibfield  {title} {\enquote {\bibinfo {title} {Switching to
  nonhyperbolic cycles from codim-2 bifurcations of equilibria in ddes},}\
  }\href@noop {} {\bibfield  {journal} {\bibinfo  {journal} {preprint}\ }
  (\bibinfo {year} {2017})}\BibitemShut {NoStop}%
\bibitem [{\citenamefont {Humphries}, \citenamefont {Calleja},\ and\
  \citenamefont {Krauskopf}(2016)}]{HCK16}%
  \BibitemOpen
  \bibfield  {author} {\bibinfo {author} {\bibfnamefont {A.}~\bibnamefont
  {Humphries}}, \bibinfo {author} {\bibfnamefont {R.}~\bibnamefont {Calleja}},
  \ and\ \bibinfo {author} {\bibfnamefont {B.}~\bibnamefont {Krauskopf}},\
  }\bibfield  {title} {\enquote {\bibinfo {title} {Resonance phenomena in a
  scalar delay differential equation with two state-dependent delays},}\
  }\href@noop {} {\bibfield  {journal} {\bibinfo  {journal} {arxiv:1607.02683}\
  } (\bibinfo {year} {2016})}\BibitemShut {NoStop}%
\bibitem [{\citenamefont {Hale}\ and\ \citenamefont {{Verduyn
  Lunel}}(1993)}]{HL93}%
  \BibitemOpen
  \bibfield  {author} {\bibinfo {author} {\bibfnamefont {J.}~\bibnamefont
  {Hale}}\ and\ \bibinfo {author} {\bibfnamefont {S.}~\bibnamefont {{Verduyn
  Lunel}}},\ }\href@noop {} {\emph {\bibinfo {title} {{Introduction to
  functional-differential equations}}}},\ \bibinfo {series} {Applied
  Mathematical Sciences}, Vol.~\bibinfo {volume} {99}\ (\bibinfo  {publisher}
  {Springer-Verlag},\ \bibinfo {address} {New York},\ \bibinfo {year} {1993})\
  pp.\ \bibinfo {pages} {x+447}\BibitemShut {NoStop}%
\bibitem [{\citenamefont {Diekmann}\ \emph {et~al.}(1995)\citenamefont
  {Diekmann}, \citenamefont {van Gils}, \citenamefont {{Verduyn Lunel}},\ and\
  \citenamefont {Walther}}]{DGLW95}%
  \BibitemOpen
  \bibfield  {author} {\bibinfo {author} {\bibfnamefont {O.}~\bibnamefont
  {Diekmann}}, \bibinfo {author} {\bibfnamefont {S.}~\bibnamefont {van Gils}},
  \bibinfo {author} {\bibfnamefont {S.}~\bibnamefont {{Verduyn Lunel}}}, \ and\
  \bibinfo {author} {\bibfnamefont {H.-O.}\ \bibnamefont {Walther}},\
  }\href@noop {} {\emph {\bibinfo {title} {{Delay equations}}}},\ \bibinfo
  {series} {Applied Mathematical Sciences}, Vol.\ \bibinfo {volume} {110}\
  (\bibinfo  {publisher} {Springer-Verlag},\ \bibinfo {address} {New York},\
  \bibinfo {year} {1995})\ pp.\ \bibinfo {pages} {xii+534}\BibitemShut
  {NoStop}%
\bibitem [{\citenamefont {Hartung}(2011)}]{Hartung11}%
  \BibitemOpen
  \bibfield  {author} {\bibinfo {author} {\bibfnamefont {F.}~\bibnamefont
  {Hartung}},\ }\bibfield  {title} {\enquote {\bibinfo {title}
  {Differentiability of solutions with respect to the initial data in
  differential equations with state-dependent delays},}\ }\href@noop {}
  {\bibfield  {journal} {\bibinfo  {journal} {J. Dyn. Diff. Eq.}\ }\textbf
  {\bibinfo {volume} {23}},\ \bibinfo {pages} {843--884} (\bibinfo {year}
  {2011})}\BibitemShut {NoStop}%
\bibitem [{\citenamefont {Walther}(2003)}]{W03}%
  \BibitemOpen
  \bibfield  {author} {\bibinfo {author} {\bibfnamefont {H.-O.}\ \bibnamefont
  {Walther}},\ }\bibfield  {title} {\enquote {\bibinfo {title} {The solution
  manifold and {$C^1$}-smoothness for differential equations with
  state-dependent delay},}\ }\href@noop {} {\bibfield  {journal} {\bibinfo
  {journal} {Journal of Differential Equations}\ }\textbf {\bibinfo {volume}
  {195}},\ \bibinfo {pages} {46--65} (\bibinfo {year} {2003})}\BibitemShut
  {NoStop}%
\bibitem [{\citenamefont {Hartung}\ \emph {et~al.}(2006)\citenamefont
  {Hartung}, \citenamefont {Krisztin}, \citenamefont {Walther},\ and\
  \citenamefont {Wu}}]{HKWW06}%
  \BibitemOpen
  \bibfield  {author} {\bibinfo {author} {\bibfnamefont {F.}~\bibnamefont
  {Hartung}}, \bibinfo {author} {\bibfnamefont {T.}~\bibnamefont {Krisztin}},
  \bibinfo {author} {\bibfnamefont {H.-O.}\ \bibnamefont {Walther}}, \ and\
  \bibinfo {author} {\bibfnamefont {J.}~\bibnamefont {Wu}},\ }\bibfield
  {title} {\enquote {\bibinfo {title} {Functional differential equations with
  state-dependent delays{\rm:} theory and applications},}\ }in\ \href@noop {}
  {\emph {\bibinfo {booktitle} {Handbook of Differential Equations{\rm:}
  Ordinary Differential Equations}}},\ Vol.~\bibinfo {volume} {3},\ \bibinfo
  {editor} {edited by\ \bibinfo {editor} {\bibfnamefont {P.}~\bibnamefont
  {Dr{\'a}bek}}, \bibinfo {editor} {\bibfnamefont {A.}~\bibnamefont
  {Ca{\~n}ada}}, \ and\ \bibinfo {editor} {\bibfnamefont {A.}~\bibnamefont
  {Fonda}}}\ (\bibinfo  {publisher} {North-Holland},\ \bibinfo {year} {2006})\
  Chap.~\bibinfo {chapter} {5}, pp.\ \bibinfo {pages} {435--545}\BibitemShut
  {NoStop}%
\bibitem [{\citenamefont {Walther}(2002)}]{W02}%
  \BibitemOpen
  \bibfield  {author} {\bibinfo {author} {\bibfnamefont {H.-O.}\ \bibnamefont
  {Walther}},\ }\bibfield  {title} {\enquote {\bibinfo {title} {Stable periodic
  motion of a system with state-dependent delay},}\ }\href@noop {} {\bibfield
  {journal} {\bibinfo  {journal} {Differential and Integral Equations}\
  }\textbf {\bibinfo {volume} {15}},\ \bibinfo {pages} {923--944} (\bibinfo
  {year} {2002})}\BibitemShut {NoStop}%
\bibitem [{\citenamefont {Insperger}, \citenamefont {Barton},\ and\
  \citenamefont {St{\'e}p{\'a}n}(2008)}]{IBS08}%
  \BibitemOpen
  \bibfield  {author} {\bibinfo {author} {\bibfnamefont {T.}~\bibnamefont
  {Insperger}}, \bibinfo {author} {\bibfnamefont {D.~A.~W.}\ \bibnamefont
  {Barton}}, \ and\ \bibinfo {author} {\bibfnamefont {G.}~\bibnamefont
  {St{\'e}p{\'a}n}},\ }\bibfield  {title} {\enquote {\bibinfo {title}
  {Criticality of {H}opf bifurcation in state-dependent delay model of turning
  processes},}\ }\href@noop {} {\bibfield  {journal} {\bibinfo  {journal}
  {International Journal of Non-Linear Mechanics}\ }\textbf {\bibinfo {volume}
  {43}},\ \bibinfo {pages} {140 -- 149} (\bibinfo {year} {2008})}\BibitemShut
  {NoStop}%
\bibitem [{\citenamefont {Insperger}, \citenamefont {St\'ep\'an},\ and\
  \citenamefont {Turi}(2007)}]{IST07}%
  \BibitemOpen
  \bibfield  {author} {\bibinfo {author} {\bibfnamefont {T.}~\bibnamefont
  {Insperger}}, \bibinfo {author} {\bibfnamefont {G.}~\bibnamefont
  {St\'ep\'an}}, \ and\ \bibinfo {author} {\bibfnamefont {J.}~\bibnamefont
  {Turi}},\ }\bibfield  {title} {\enquote {\bibinfo {title} {State-dependent
  delay in regenerative turning processes},}\ }\href@noop {} {\bibfield
  {journal} {\bibinfo  {journal} {Nonlinear Dynamics}\ }\textbf {\bibinfo
  {volume} {47}},\ \bibinfo {pages} {275--283} (\bibinfo {year}
  {2007})}\BibitemShut {NoStop}%
\bibitem [{\citenamefont {Luca}\ \emph {et~al.}(2010)\citenamefont {Luca},
  \citenamefont {Guglielmi}, \citenamefont {Humphries},\ and\ \citenamefont
  {Politi}}]{DGHP10}%
  \BibitemOpen
  \bibfield  {author} {\bibinfo {author} {\bibfnamefont {J.~D.}\ \bibnamefont
  {Luca}}, \bibinfo {author} {\bibfnamefont {N.}~\bibnamefont {Guglielmi}},
  \bibinfo {author} {\bibfnamefont {A.}~\bibnamefont {Humphries}}, \ and\
  \bibinfo {author} {\bibfnamefont {A.}~\bibnamefont {Politi}},\ }\bibfield
  {title} {\enquote {\bibinfo {title} {Electromagnetic two-body problem{\rm:}
  recurrent dynamics in the presence of state-dependent delay},}\ }\href@noop
  {} {\bibfield  {journal} {\bibinfo  {journal} {Journal of Physics A}\
  }\textbf {\bibinfo {volume} {43}} (\bibinfo {year} {2010})}\BibitemShut
  {NoStop}%
\bibitem [{\citenamefont {Pyragas}\ and\ \citenamefont {Pyragas}(2011)}]{PP11}%
  \BibitemOpen
  \bibfield  {author} {\bibinfo {author} {\bibfnamefont {V.}~\bibnamefont
  {Pyragas}}\ and\ \bibinfo {author} {\bibfnamefont {K.}~\bibnamefont
  {Pyragas}},\ }\bibfield  {title} {\enquote {\bibinfo {title} {Adaptive
  modification of the delayed feedback control algorithm with a continuously
  varying time delay},}\ }\href@noop {} {\bibfield  {journal} {\bibinfo
  {journal} {Physics Letters A}\ }\textbf {\bibinfo {volume} {375}},\ \bibinfo
  {pages} {3866--3871} (\bibinfo {year} {2011})}\BibitemShut {NoStop}%
\bibitem [{\citenamefont {Craig}, \citenamefont {Humphries},\ and\
  \citenamefont {M.C.Mackey}(2016)}]{CHM16}%
  \BibitemOpen
  \bibfield  {author} {\bibinfo {author} {\bibfnamefont {M.}~\bibnamefont
  {Craig}}, \bibinfo {author} {\bibfnamefont {A.}~\bibnamefont {Humphries}}, \
  and\ \bibinfo {author} {\bibnamefont {M.C.Mackey}},\ }\bibfield  {title}
  {\enquote {\bibinfo {title} {A mathematical model of granulopoiesis
  incorporating the negative feedback dynamics and kinetics of
  {G-CSF}/neutrophil binding and internalisation},}\ }\href@noop {} {\bibfield
  {journal} {\bibinfo  {journal} {Bulletin of Mathematical Biology}\ }\textbf
  {\bibinfo {volume} {78}},\ \bibinfo {pages} {2304--2357} (\bibinfo {year}
  {2016})}\BibitemShut {NoStop}%
\bibitem [{\citenamefont {Guglielmi}\ and\ \citenamefont
  {Hairer}(2001)}]{GH01}%
  \BibitemOpen
  \bibfield  {author} {\bibinfo {author} {\bibfnamefont {N.}~\bibnamefont
  {Guglielmi}}\ and\ \bibinfo {author} {\bibfnamefont {E.}~\bibnamefont
  {Hairer}},\ }\bibfield  {title} {\enquote {\bibinfo {title} {Implementing
  {R}adau {IIA} methods for stiff delay differential equations},}\ }\href@noop
  {} {\bibfield  {journal} {\bibinfo  {journal} {Computing}\ }\textbf {\bibinfo
  {volume} {67}},\ \bibinfo {pages} {1--12} (\bibinfo {year}
  {2001})}\BibitemShut {NoStop}%
\bibitem [{\citenamefont {Sieber}(2012)}]{S12}%
  \BibitemOpen
  \bibfield  {author} {\bibinfo {author} {\bibfnamefont {J.}~\bibnamefont
  {Sieber}},\ }\bibfield  {title} {\enquote {\bibinfo {title} {Finding periodic
  orbits in state-dependent delay differential equations as roots of algebraic
  equations},}\ }\href@noop {} {\bibfield  {journal} {\bibinfo  {journal}
  {Discrete and Continuous Dynamical Systems A}\ }\textbf {\bibinfo {volume}
  {32}},\ \bibinfo {pages} {2607--2651} (\bibinfo {year} {2012})}\BibitemShut
  {NoStop}%
\bibitem [{\citenamefont {Engelborghs}, \citenamefont {Luzyanina},\ and\
  \citenamefont {Roose}(2002)}]{ELR02}%
  \BibitemOpen
  \bibfield  {author} {\bibinfo {author} {\bibfnamefont {K.}~\bibnamefont
  {Engelborghs}}, \bibinfo {author} {\bibfnamefont {T.}~\bibnamefont
  {Luzyanina}}, \ and\ \bibinfo {author} {\bibfnamefont {D.}~\bibnamefont
  {Roose}},\ }\bibfield  {title} {\enquote {\bibinfo {title} {Numerical
  bifurcation analysis of delay differential equations using {DDE-BIFTOOL}},}\
  }\href@noop {} {\bibfield  {journal} {\bibinfo  {journal} {ACM Transactions
  on Mathematical Software}\ }\textbf {\bibinfo {volume} {28}},\ \bibinfo
  {pages} {1--21} (\bibinfo {year} {2002})}\BibitemShut {NoStop}%
\bibitem [{\citenamefont {Engelborghs}, \citenamefont {Luzyanina},\ and\
  \citenamefont {Samaey}(2001)}]{ELS01}%
  \BibitemOpen
  \bibfield  {author} {\bibinfo {author} {\bibfnamefont {K.}~\bibnamefont
  {Engelborghs}}, \bibinfo {author} {\bibfnamefont {T.}~\bibnamefont
  {Luzyanina}}, \ and\ \bibinfo {author} {\bibfnamefont {G.}~\bibnamefont
  {Samaey}},\ }\href@noop {} {\enquote {\bibinfo {title} {{{DDE-BIFTOOL}
  v.2.00{\rm:} a {Matlab} package for bifurcation analysis of delay
  differential equations}},}\ }\bibinfo {type} {Report TW}\ \bibinfo {number}
  {330}\ (\bibinfo  {institution} {Katholieke Universiteit Leuven},\ \bibinfo
  {year} {2001})\BibitemShut {NoStop}%
\bibitem [{\citenamefont {Sieber}\ \emph {et~al.}()\citenamefont {Sieber},
  \citenamefont {Engelborghs}, \citenamefont {Luzyanina}, \citenamefont
  {Samaey},\ and\ \citenamefont {Roose}}]{ddebiftoolmanual}%
  \BibitemOpen
  \bibfield  {author} {\bibinfo {author} {\bibfnamefont {J.}~\bibnamefont
  {Sieber}}, \bibinfo {author} {\bibfnamefont {K.}~\bibnamefont {Engelborghs}},
  \bibinfo {author} {\bibfnamefont {T.}~\bibnamefont {Luzyanina}}, \bibinfo
  {author} {\bibfnamefont {G.}~\bibnamefont {Samaey}}, \ and\ \bibinfo {author}
  {\bibfnamefont {D.}~\bibnamefont {Roose}},\ }\href@noop {} {\emph {\bibinfo
  {title} {DDE-BIFTOOL Manual --- Bifurcation analysis of delay differential
  equations}}},\ \bibinfo {address}
  {\url{sourceforge.net/projects/ddebiftool}}\BibitemShut {NoStop}%
\bibitem [{\citenamefont {Thompson}\ and\ \citenamefont
  {Stewart}(2002)}]{TS02}%
  \BibitemOpen
  \bibfield  {author} {\bibinfo {author} {\bibfnamefont {J.~M.~T.}\
  \bibnamefont {Thompson}}\ and\ \bibinfo {author} {\bibfnamefont {H.~B.}\
  \bibnamefont {Stewart}},\ }\href@noop {} {\emph {\bibinfo {title} {{Nonlinear
  dynamics and chaos}}}},\ \bibinfo {edition} {2nd}\ ed.\ (\bibinfo
  {publisher} {Wiley},\ \bibinfo {address} {Chichester, UK},\ \bibinfo {year}
  {2002})\BibitemShut {NoStop}%
\bibitem [{\citenamefont {Guckenheimer}\ and\ \citenamefont
  {Holmes}(1990)}]{GH83}%
  \BibitemOpen
  \bibfield  {author} {\bibinfo {author} {\bibfnamefont {J.}~\bibnamefont
  {Guckenheimer}}\ and\ \bibinfo {author} {\bibfnamefont {P.}~\bibnamefont
  {Holmes}},\ }\href@noop {} {\emph {\bibinfo {title} {{Nonlinear oscillations,
  dynamical systems, and bifurcations of vector fields}}}},\ \bibinfo {series}
  {Applied Mathematical Sciences}, Vol.~\bibinfo {volume} {42}\ (\bibinfo
  {publisher} {Springer-Verlag},\ \bibinfo {address} {New York},\ \bibinfo
  {year} {1990})\ pp.\ \bibinfo {pages} {xvi+459}\BibitemShut {NoStop}%
\bibitem [{\citenamefont {Kuznetsov}(2004)}]{K04}%
  \BibitemOpen
  \bibfield  {author} {\bibinfo {author} {\bibfnamefont {Y.~A.}\ \bibnamefont
  {Kuznetsov}},\ }\href@noop {} {\emph {\bibinfo {title} {{Elements of Applied
  Bifurcation Theory}}}},\ \bibinfo {edition} {3rd}\ ed.,\ \bibinfo {series}
  {Applied Mathematical Sciences}, Vol.\ \bibinfo {volume} {112}\ (\bibinfo
  {publisher} {Springer-Verlag},\ \bibinfo {address} {New York},\ \bibinfo
  {year} {2004})\ pp.\ \bibinfo {pages} {xxii+631}\BibitemShut {NoStop}%
\bibitem [{\citenamefont {Stumpf}(2011)}]{Stu11}%
  \BibitemOpen
  \bibfield  {author} {\bibinfo {author} {\bibfnamefont {E.}~\bibnamefont
  {Stumpf}},\ }\bibfield  {title} {\enquote {\bibinfo {title} {The existence
  and $c^{1}$-smoothness of local center-unstable manifolds for differential
  equations with state-dependent delay},}\ }\href@noop {} {\bibfield  {journal}
  {\bibinfo  {journal} {Rostocker Mathematisches Kolloquium}\ }\textbf
  {\bibinfo {volume} {66}},\ \bibinfo {pages} {3--44} (\bibinfo {year}
  {2011})}\BibitemShut {NoStop}%
\bibitem [{\citenamefont {Stumpf}(2015{\natexlab{a}})}]{Stu15a}%
  \BibitemOpen
  \bibfield  {author} {\bibinfo {author} {\bibfnamefont {E.}~\bibnamefont
  {Stumpf}},\ }\bibfield  {title} {\enquote {\bibinfo {title} {Attraction
  property of local center-unstable manifolds for differential equations with
  state-dependent delay},}\ }\href@noop {} {\bibfield  {journal} {\bibinfo
  {journal} {Electronic Journal of Qualitative Theory of Differential
  Equations}\ }\textbf {\bibinfo {volume} {2015}},\ \bibinfo {pages} {1--45}
  (\bibinfo {year} {2015}{\natexlab{a}})}\BibitemShut {NoStop}%
\bibitem [{\citenamefont {Walther}(2004)}]{W04}%
  \BibitemOpen
  \bibfield  {author} {\bibinfo {author} {\bibfnamefont {H.-O.}\ \bibnamefont
  {Walther}},\ }\bibfield  {title} {\enquote {\bibinfo {title} {Smoothness
  properties of semiflows for differential equations with state-dependent
  delays},}\ }\href@noop {} {\bibfield  {journal} {\bibinfo  {journal} {Journal
  of Mathematical Sciences}\ }\textbf {\bibinfo {volume} {124}},\ \bibinfo
  {pages} {5193--5207} (\bibinfo {year} {2004})}\BibitemShut {NoStop}%
\bibitem [{\citenamefont {Stumpf}(2015{\natexlab{b}})}]{Stu15}%
  \BibitemOpen
  \bibfield  {author} {\bibinfo {author} {\bibfnamefont {E.}~\bibnamefont
  {Stumpf}},\ }\bibfield  {title} {\enquote {\bibinfo {title} {Local stability
  analysis of differential equations with state-dependent delay},}\ }\href@noop
  {} {\bibfield  {journal} {\bibinfo  {journal} {arXiv:1502.03142}\ } (\bibinfo
  {year} {2015}{\natexlab{b}})}\BibitemShut {NoStop}%
\bibitem [{\citenamefont {Krisztin}(2003)}]{K03}%
  \BibitemOpen
  \bibfield  {author} {\bibinfo {author} {\bibfnamefont {T.}~\bibnamefont
  {Krisztin}},\ }\bibfield  {title} {\enquote {\bibinfo {title} {A local
  unstable manifold for differential equations with state-dependent delay},}\
  }\href@noop {} {\bibfield  {journal} {\bibinfo  {journal} {Discrete Contin.
  Dynam. Systems}\ }\textbf {\bibinfo {volume} {9}},\ \bibinfo {pages}
  {993--1028} (\bibinfo {year} {2003})}\BibitemShut {NoStop}%
\bibitem [{\citenamefont {Stumpf}\ \emph {et~al.}(2016)\citenamefont {Stumpf}
  \emph {et~al.}}]{Stu16b}%
  \BibitemOpen
  \bibfield  {author} {\bibinfo {author} {\bibfnamefont {E.}~\bibnamefont
  {Stumpf}} \emph {et~al.},\ }\bibfield  {title} {\enquote {\bibinfo {title} {A
  note on local center manifolds for differential equations with
  state-dependent delay},}\ }\href@noop {} {\bibfield  {journal} {\bibinfo
  {journal} {Differential and Integral Equations}\ }\textbf {\bibinfo {volume}
  {29}},\ \bibinfo {pages} {1093--1106} (\bibinfo {year} {2016})}\BibitemShut
  {NoStop}%
\bibitem [{\citenamefont {Krisztin}(2006{\natexlab{a}})}]{K06}%
  \BibitemOpen
  \bibfield  {author} {\bibinfo {author} {\bibfnamefont {T.}~\bibnamefont
  {Krisztin}},\ }\bibfield  {title} {\enquote {\bibinfo {title} {Smooth center
  manifolds for differential equations with state-dependent delay},}\ }in\
  \href@noop {} {\emph {\bibinfo {booktitle} {AIMS Conference Poitiers}}}\
  (\bibinfo {year} {2006})\BibitemShut {NoStop}%
\bibitem [{\citenamefont {Krisztin}(2006{\natexlab{b}})}]{K06a}%
  \BibitemOpen
  \bibfield  {author} {\bibinfo {author} {\bibfnamefont {T.}~\bibnamefont
  {Krisztin}},\ }\bibfield  {title} {\enquote {\bibinfo {title}
  {{$C^1$}-smoothness of center manifolds for differential equations with
  state-dependent delay},}\ }in\ \href@noop {} {\emph {\bibinfo {booktitle}
  {Nonlinear dynamics and evolution equations}}},\ Vol.~\bibinfo {volume} {48}\
  (\bibinfo  {publisher} {Fields Institute Communications},\ \bibinfo {year}
  {2006})\ pp.\ \bibinfo {pages} {213--226}\BibitemShut {NoStop}%
\bibitem [{\citenamefont {Eichmann}(2006)}]{E06}%
  \BibitemOpen
  \bibfield  {author} {\bibinfo {author} {\bibfnamefont {M.}~\bibnamefont
  {Eichmann}},\ }\emph {\bibinfo {title} {A local {H}opf Bifurcation Theorem
  for differential equations with state-dependent delays}},\ \href@noop {}
  {Ph.D. thesis},\ \bibinfo  {school} {University of Giessen} (\bibinfo {year}
  {2006})\BibitemShut {NoStop}%
\bibitem [{\citenamefont {Sieber}(2013)}]{S13}%
  \BibitemOpen
  \bibfield  {author} {\bibinfo {author} {\bibfnamefont {J.}~\bibnamefont
  {Sieber}},\ }\bibfield  {title} {\enquote {\bibinfo {title} {Extended systems
  for delay-differential equations as implemented in the extensions to
  {DDE}-{B}if{T}ool},}\ }\href@noop {} {\bibfield  {journal} {\bibinfo
  {journal} {Figshare}\ }\textbf {\bibinfo {volume}
  {\url{http://dx.doi.org/10.6084/m9.figshare.757725}}} (\bibinfo {year}
  {2013})}\BibitemShut {NoStop}%
\bibitem [{\citenamefont {Tziperman}\ \emph {et~al.}(1998)\citenamefont
  {Tziperman}, \citenamefont {Cane}, \citenamefont {Zebiak}, \citenamefont
  {Xue},\ and\ \citenamefont {Blumenthal}}]{TCZXB98}%
  \BibitemOpen
  \bibfield  {author} {\bibinfo {author} {\bibfnamefont {E.}~\bibnamefont
  {Tziperman}}, \bibinfo {author} {\bibfnamefont {M.}~\bibnamefont {Cane}},
  \bibinfo {author} {\bibfnamefont {S.}~\bibnamefont {Zebiak}}, \bibinfo
  {author} {\bibfnamefont {Y.}~\bibnamefont {Xue}}, \ and\ \bibinfo {author}
  {\bibfnamefont {B.}~\bibnamefont {Blumenthal}},\ }\bibfield  {title}
  {\enquote {\bibinfo {title} {Locking of {E}l {N}i{\~n}o's peak time to the
  end of the calendar year in the delayed oscillator picture of {ENSO}},}\
  }\href@noop {} {\bibfield  {journal} {\bibinfo  {journal} {Journal of
  Climate}\ }\textbf {\bibinfo {volume} {11}},\ \bibinfo {pages} {2191--2199}
  (\bibinfo {year} {1998})}\BibitemShut {NoStop}%
\bibitem [{\citenamefont {Keane}, \citenamefont {Krauskopf},\ and\
  \citenamefont {Postlethwaite}(2015)}]{KKP15}%
  \BibitemOpen
  \bibfield  {author} {\bibinfo {author} {\bibfnamefont {A.}~\bibnamefont
  {Keane}}, \bibinfo {author} {\bibfnamefont {B.}~\bibnamefont {Krauskopf}}, \
  and\ \bibinfo {author} {\bibfnamefont {C.}~\bibnamefont {Postlethwaite}},\
  }\bibfield  {title} {\enquote {\bibinfo {title} {Delayed feedback versus
  seasonal forcing: Resonance phenomena in an {E}l {N}i{\~n}o {S}outhern
  {O}scillation model},}\ }\href@noop {} {\bibfield  {journal} {\bibinfo
  {journal} {SIAM Journal on Applied Dynamical Systems}\ }\textbf {\bibinfo
  {volume} {14}},\ \bibinfo {pages} {1229--1257} (\bibinfo {year}
  {2015})}\BibitemShut {NoStop}%
\bibitem [{\citenamefont {Keane}, \citenamefont {Krauskopf},\ and\
  \citenamefont {Postlethwaite}(2016)}]{KKP16}%
  \BibitemOpen
  \bibfield  {author} {\bibinfo {author} {\bibfnamefont {A.}~\bibnamefont
  {Keane}}, \bibinfo {author} {\bibfnamefont {B.}~\bibnamefont {Krauskopf}}, \
  and\ \bibinfo {author} {\bibfnamefont {C.}~\bibnamefont {Postlethwaite}},\
  }\bibfield  {title} {\enquote {\bibinfo {title} {Investigating irregular
  behavior in a model for the {E}l {N}i{\~n}o {S}outhern {O}scillation with
  positive and negative delayed feedback},}\ }\href@noop {} {\bibfield
  {journal} {\bibinfo  {journal} {SIAM Journal on Applied Dynamical Systems}\
  }\textbf {\bibinfo {volume} {15}},\ \bibinfo {pages} {1656--1689} (\bibinfo
  {year} {2016})}\BibitemShut {NoStop}%
\bibitem [{\citenamefont {Krauskopf}\ and\ \citenamefont
  {Sieber}(2014)}]{KS14}%
  \BibitemOpen
  \bibfield  {author} {\bibinfo {author} {\bibfnamefont {B.}~\bibnamefont
  {Krauskopf}}\ and\ \bibinfo {author} {\bibfnamefont {J.}~\bibnamefont
  {Sieber}},\ }\bibfield  {title} {\enquote {\bibinfo {title} {Bifurcation
  analysis of delay-induced resonances of the {E}l {N}i{\~n}o {S}outhern
  {O}scillation},}\ }\href@noop {} {\bibfield  {journal} {\bibinfo  {journal}
  {Proc. Roy. Soc. London A}\ }\textbf {\bibinfo {volume} {470}} (\bibinfo
  {year} {2014})}\BibitemShut {NoStop}%
\bibitem [{\citenamefont {Kuznetsov}(1999)}]{K99}%
  \BibitemOpen
  \bibfield  {author} {\bibinfo {author} {\bibfnamefont {Y.~A.}\ \bibnamefont
  {Kuznetsov}},\ }\bibfield  {title} {\enquote {\bibinfo {title} {Numerical
  normalization techniques for all codim 2 bifurcations of equilibria in
  ode's},}\ }\href@noop {} {\bibfield  {journal} {\bibinfo  {journal} {SIAM
  journal on numerical analysis}\ }\textbf {\bibinfo {volume} {36}},\ \bibinfo
  {pages} {1104--1124} (\bibinfo {year} {1999})}\BibitemShut {NoStop}%
\bibitem [{\citenamefont {van Gils}\ \emph {et~al.}(2013)\citenamefont {van
  Gils}, \citenamefont {Janssens}, \citenamefont {Kuznetsov},\ and\
  \citenamefont {Visser}}]{GJKV13}%
  \BibitemOpen
  \bibfield  {author} {\bibinfo {author} {\bibfnamefont {S.~A.}\ \bibnamefont
  {van Gils}}, \bibinfo {author} {\bibfnamefont {S.~G.}\ \bibnamefont
  {Janssens}}, \bibinfo {author} {\bibfnamefont {Y.~A.}\ \bibnamefont
  {Kuznetsov}}, \ and\ \bibinfo {author} {\bibfnamefont {S.}~\bibnamefont
  {Visser}},\ }\bibfield  {title} {\enquote {\bibinfo {title} {On local
  bifurcations in neural field models with transmission delays},}\ }\href@noop
  {} {\bibfield  {journal} {\bibinfo  {journal} {Journal of mathematical
  biology}\ }\textbf {\bibinfo {volume} {66}},\ \bibinfo {pages} {837--887}
  (\bibinfo {year} {2013})}\BibitemShut {NoStop}%
\bibitem [{\citenamefont {Hartung}\ and\ \citenamefont
  {Turi}(2001)}]{hartung2001linearized}%
  \BibitemOpen
  \bibfield  {author} {\bibinfo {author} {\bibfnamefont {F.}~\bibnamefont
  {Hartung}}\ and\ \bibinfo {author} {\bibfnamefont {J.}~\bibnamefont {Turi}},\
  }\bibfield  {title} {\enquote {\bibinfo {title} {Linearized stability in
  functional differential equations with state-dependent delays},}\ }in\
  \href@noop {} {\emph {\bibinfo {booktitle} {2000 International Conference on
  Dynamical Systems and Differential Equations}}}\ (\bibinfo {organization}
  {American Institute of Mathematical Sciences},\ \bibinfo {year}
  {2001})\BibitemShut {NoStop}%
\bibitem [{\citenamefont {Cooke}\ and\ \citenamefont
  {Huang}(1996)}]{cooke1996problem}%
  \BibitemOpen
  \bibfield  {author} {\bibinfo {author} {\bibfnamefont {K.~L.}\ \bibnamefont
  {Cooke}}\ and\ \bibinfo {author} {\bibfnamefont {W.}~\bibnamefont {Huang}},\
  }\bibfield  {title} {\enquote {\bibinfo {title} {On the problem of
  linearization for state-dependent delay differential equations},}\
  }\href@noop {} {\bibfield  {journal} {\bibinfo  {journal} {Proceedings of the
  American Mathematical Society}\ ,\ \bibinfo {pages} {1417--1426}} (\bibinfo
  {year} {1996})}\BibitemShut {NoStop}%
\end{thebibliography}%

\end{document}